\documentclass[a4paper,12pt]{article}
\usepackage[utf8]{inputenc}
\usepackage[esperanto,portuguese,english]{babel}
\usepackage{epsfig}
\usepackage{icomma}
\usepackage{float}

\newcommand{\ppparallel}[1]{} \newcommand{\ppmaldekstra}[1]{} \newcommand{\ppdekstra}[1]{}
\usepackage{parallel} \renewcommand{\ppparallel}[1]{#1} 

\newcommand{\comentario}[1]{}

\newcommand{\dd}{\mathrm{d}}

\ppparallel{
\newlength{\pplw}\setlength{\pplw}{0.47\textwidth}
\newlength{\pprw}\setlength{\pprw}{0.49\textwidth}

\newcommand{\ppp}{\ParallelPar}
\newcommand{\ppn}{\noindent}              
\newcommand{\ppl}[1]{\ParallelLText{\selectlanguage{esperanto}#1}}
\newcommand{\ppr}[1]{\ParallelRText{\selectlanguage{english}#1}\ppp}
\newcommand{\ppln}[1]
{\ParallelLText{\ppn \selectlanguage{esperanto}#1}} 
\newcommand{\pprn}[1]
{\ParallelRText{\ppn \selectlanguage{english}#1}\ppp} 

\newcommand{\ppsection}[3][0ex]{\vspace{2em}
\ppl{\section{#2} \vspace{#1}} \ppa \nopagebreak
\ppR{\section{#3}} \ppp \nopagebreak}

\newcommand{\bea}{\vspace{-1ex}\begin{eqnarray}}
\newcommand{\eea}{\end{eqnarray}}
}

\ppmaldekstra{
\newcommand{\ppl}[1]{\selectlanguage{esperanto}#1}
\newcommand{\ppln}[1]{\noindent \selectlanguage{esperanto}#1}
\newcommand{\ppr}[1]{\selectlanguage{english}}
\newcommand{\pprn}[1]{\noindent \selectlanguage{english}}
\newcommand{\ppsection}[3][0ex]{\section{#2}}

\newcommand{\bea}{\begin{eqnarray}}
\newcommand{\eea}{\end{eqnarray}}
}

\ppdekstra{
\newcommand{\ppl}[1]{\selectlanguage{esperanto}}
\newcommand{\ppln}[1]{\noindent \selectlanguage{esperanto}}
\newcommand{\ppr}[1]{\selectlanguage{english}#1}
\newcommand{\pprn}[1]{\noindent \selectlanguage{english}#1}
\newcommand{\ppsection}[3][0ex]{\section{#3}}

\newcommand{\bea}{\begin{eqnarray}}
\newcommand{\eea}{\end{eqnarray}}
}

\title{{\bf Schwarz-Christoffel: piliero en rivero \ppparallel{\\ Schwarz-Christoffel: a pillar in a river}}}
\author{F.M. Paiva \\
{\small Departamento de F\'isica, Campus Humait\'a II, Col\'egio Pedro II} \\
{\small Rua Humait\'a 80, 22261-001  Rio de Janeiro-RJ, Brasil; fmpaiva@cbpf.br}
\vspace{.7ex} \\
A.F.F. Teixeira \\
{\small Centro Brasileiro de Pesquisas F\'isicas} \\
{\small 22290-180 Rio de Janeiro-RJ, Brasil; teixeira@cbpf.br}}
\date{}

\begin{document}
\selectlanguage{esperanto}
\maketitle
\thispagestyle{empty}

\begin{abstract}\selectlanguage{esperanto}
La transformoj de Schwarz-Christoffel mapas, konforme, la superan kompleksan duon-ebenon al regiono limigita per rektaj segmentoj. \^Ci tie ni priskribas kiel konvene kunigi mapon de la suba duon-ebeno al mapo de la supera duon-ebeno. Ni emfazas la bezonon de klara difino de angulo de kompleksa nombro, por tiu kunigo. Ni diskutas kelkajn ekzemplojn kaj donas interesan aplikon pri movado de fluido. 

\ppparallel{\selectlanguage{english}
Schwarz-Christoffel transformations map, conformally, the complex upper half plane into a region bounded by right segments. Here we describe how to couple conveniently a map of the lower half plane to the map of the upper half plane. We emphasize the need of a clear definition of angle of a complex, to that coupling. We discuss some examples and give an interesting application for motion of fluid.} 

\end{abstract}

\ppparallel{
\begin{Parallel}[v]{\pplw}{\pprw}
}

\ppparallel{\section*{\vspace{-2em}}\vspace{-2ex}}   

\ppsection[0.6ex]{\label{intro}Enkonduko}{Introduction} 

\ppln{Ni konsideru mapojn de kartezia ebeno $z$, a\u u de {\it regiono} el \^gi. Ili asocios, al \^ciu punkto $z=[x,y]$ de la regiono, unu {\it imagan} punkton $w=[u,v]$ en kartezia ebeno $w$.}
\pprn{Let's consider maps of the cartesian plane $z$, or of {\it a region} of it. They will associate, to each point $z=[x,y]$ of the region, an {\it image} point $w=[u,v]$ in a cartesian plane $w$.}

\ppl{Niaj mapoj estos tre specialaj: ili estos {\it konfomaj}. Per tiuj mapoj, du linioj en ebeno $z$ kiuj krucas je angulo $\alpha$ havos imagojn krucantajn je la sama $\alpha$. Sekvas, ke tiuj transformoj mapas malgrandan geometrian bildon en alia malgranda geometria bildo kun sama formo, tamen \^generale kun alia orienti\^go kaj alia grandeco \mbox{\cite[pa\^go 541]{Hildebrand}}. Figuro~\ref{160315a} montras mapon konforman kaj alian nekonforman.}
\ppr{Our maps will be very special: they will be {\em conformal}. In these maps, two lines in the plane $z$ crossing with angle $\alpha$ will have images also crossing with $\alpha$. It follows that these transformations map a small geome\-tric figure into another small geome\-tric figure with the same form, though generally with other orientation and other size \mbox{\cite[page ~541]{Hildebrand}}. \mbox{Figure~\ref{160315a}} shows a conformal map, and another non-conformal.}

\begin{figure}[H]                                                                         
\centerline{\epsfig{file=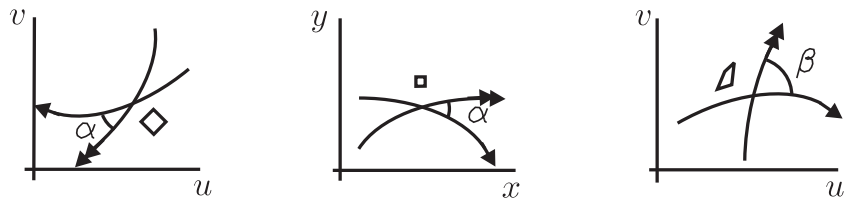,width=11cm}} 
\caption{Centre, kruca\^\j o de du linioj en ebeno $[x,y]$, kaj malgranda kvadrato; maldekstre, iliaj imagoj per iu konforma mapo; dekstre, iliaj imagoj per iu ne-konforma mapo, estante $\beta\neq\alpha$.
\newline
Figure~\ref{160315a}: In the center, crossing of two lines in the plane $[x,y]$, and a small square; on the left, their images under a conformal map; on the right, their images under some non-conformal map, being $\beta\neq\alpha$.} 
\label{160315a} 
\end{figure}

\ppl{Je ajn mapo (konforma a\u u ne), la koordinatoj $u$ kaj $v$ de \^ciu imaga punkto estas funkcioj de koordinatoj $x$ kaj $y$ de la responda anta\u u-imaga punkto: $u=u(x,y)$ kaj \mbox{$v=v(x,y)$}. Demando: kiujn matematikajn propretojn tiuj funkcioj $u$ kaj $v$ havas, tial ke la mapo $[x,y]\rightarrow[u,v]$ estu konforma? Respondo: kondi\^cojn de Cauchy-Riemann \mbox{\cite[p. 46]{Knopp},}}
\ppr{Under any map (conformal or not), the coordinates $u$ and $v$ of each image point are functions of the coordinates $x$ and $y$ of the corresponding pre-image point: $u=u(x,y)$ and \mbox{$v=v(x,y)$}. Question: what mathematical properties these functions $u$ and $v$ have, so that the map $[x,y]\rightarrow[u,\,v]$ be conformal? Answer: the conditions of Cauchy-Riemann \mbox{\cite[p. 46]{Knopp},}}
\bea                                                                                  \label{CR}
\frac{\partial u}{\partial x}=\frac{\partial v}{\partial y}\,, \hspace{5mm}  \frac{\partial u}{\partial y}=-\frac{\partial v}{\partial x}\,.  
\eea  
\ppln{Ekzemple, la mapo $u(x,y)=x^2-y^2, \hspace{1mm} v(x,y)\!=2xy$ estas konforma, kontra\u ue \mbox{$u(x,y)=x+y$,} $v(x,y)=x-y$ ne estas.}
\pprn{For instance, the map $u(x,y)=x^2-y^2\,,\,\, v(x,y)=2xy$ is conformal, while \mbox{$u(x,y)=x+y$}, $v(x,y)=x-y$ is not.}

\ppl{Oportune, kondi\^coj (\ref{CR}) implicas, ke $u$ kaj $v$ amba\u u estas {\it harmoniaj} funkcioj: }
\ppr{By the way, the conditions (\ref{CR}) imply $u$ and $v$ are {\it harmonic}, both:}

\bea                                                                                \label{harm}
\Delta u\equiv\frac{\partial}{\partial x}\left(\frac{\partial u}{\partial x}\right) + \frac{\partial}{\partial y}\left(\frac{\partial u}{\partial y}\right) =0\,, 
\hspace{5mm} 
\Delta v\equiv\frac{\partial}{\partial x}\left(\frac{\partial v}{\partial x}\right) + \frac{\partial}{\partial y}\left(\frac{\partial v}{\partial y}\right) =0\,.
\eea

\ppln{Se du harmoniaj funkcioj, $u$ kaj $v$, plenumas kondi\^cojn (\ref{CR}), oni diras ke $v$ estas {\it harmonia dualo} de $u$, kaj ke $u$ estas harmonia dualo de $-v$. Harmoniaj funkcioj havas tre interesajn proprecojn, kaj ni beda\u uras ke ni ne povas halti\^gi por plezuri\^gi pri tio \mbox{\cite[p. 508]{Needham}.}}
\pprn{If two harmonic functions, $u$ and $v$, obey conditions (\ref{CR}), one says that $v$ is {\it dual harmonic} of $u$, and that $u$ is dual harmonic of $-v$. Harmonic functions have very interesting properties, and we regret to have not enough space and time to appreciate them here \mbox{\cite[p. 508]{Needham}.}}

\ppl{Ni uzos la kompleksan kalkulon, por faciligi traktadon de konformaj mapoj. Por tio, al \^ciu punkto $z=[x,y]$ de la kartezia ebeno $z$ ni asocios la kompleksan nombron $z=x+i\,y$; same, al \^ciu punkto $w=[u,v]$ de la kartezia ebeno $w$ ni asocios la kompleksan nombron $w=u+i\,v$. Plue, havante paron de ajn realaj funkcioj, $u(x,y)$ kaj $v(x,y)$, ni difinos la kompleksan funkcion $w(x,y)=u(x,y)+i\,v(x,y)$; \^ci tiu funkcio asocias, al \^ciu punkto $[x,y]$ en la ebeno $z$, unu punkton $[u,v]$ en la ebeno $w$.}
\ppr{We use the complex calculus, to facilitate handling of conformal maps. To that end, to each point $z=[x,y]$ of the cartesian plane $z$ we associate the complex number $z=x+i\,y$; equally, to each point $w=[u,v]$ of the cartesian plane $w$ we associate the complex number $w=u+i\,v$. Still, having a pair of arbitrary real functions, $u(x,y)$ and $v(x,y)$, we define the complex function $w(x,y)=u(x,y)+i\,v(x,y)$; this function associates, to each point $[x,y]$ of plane $z$, a point $[u,v]$ in plane $w$.}

\ppl{Se funkcio $u(x,y)$ estas harmonia, kaj se $v(x,y)$ estas \^gia harmonia dualo, tiuokaze la funkcio $w(x,y)=u(x,y)+i\,v(x,y)$ estos dirita {\it kompleksa konforma}. Okazas, ke se ni anstata\u uigas $x$ per $\frac{1}{2}(z+z^*)$\, kaj\, $y$ per $\frac{1}{2\,i}(z-z^*)$ en la kompleksa konforma funkcio $w(x,y)$, ni aperigos la {\it analitikan funkcion} $w(z)$, libera de $z^*$ (la konjuga\^\j on de $z$).}
\ppr{If the function $u(x,y)$ is harmonic, and if $v(x,y)$ is its dual harmonic, then the function $w(x,y)=u(x,y)+i\,v(x,y)$ is told {\it complex conformal}. It happens that, if we substitute $x\rightarrow\frac{1}{2}(z+z^*)$\, and\, $y\rightarrow\frac{1}{2\,i}(z-z^*)$ in the complex conformal function $w(x,y)$, we make appear the {\it analytical function} $w(z)$, exempt of $z^*$ (the complex conjugate of $z$).}

\ppl{Ni memoru, ke iu kompleksa nombro $\Omega$ skribi\^gas en polara formo kiel $\Omega=\rho\,{\rm e}^{i\,\alpha}$, estante $\rho$ reala ne-negativa. La angulo $\alpha$, anka\u u reala, indikas la {\it orienti\^gon} de la vektoro $0\rightarrow\Omega$ en la kompleksa ebeno. {\it En \^ci tiu teksto} ni konvencias $\alpha\in(-\pi,\,\pi\,]$ kaj skribas $\alpha\,=\,$\protect\angle\,$\Omega$.}
\ppr{We remember that a complex number $\Omega$ is written in polar form as $\Omega=\rho\,{\rm e}^{i\,\alpha}$, with $\rho$ real non-negative. The angle $\alpha$, also real, indicates the {\it orientation} of the vector $0\rightarrow\Omega$ in the complex plane. {\it In this text} we agree $\alpha\in(-\pi,\,\pi\,]$ and write $\alpha\,=\,$\protect\angle\,$\Omega$.} 

\ppl{Ni anka\u u memoru, ke se $\Omega=\Omega_1\Omega_2\,$, tio estas, $\rho\,{\rm e}^{i\,\alpha}=(\rho_1\,{\rm e}^{i\,\alpha_1})(\rho_2\,{\rm e}^{i\,\alpha_2})$, tiel $\rho=\rho_1\rho_2$ kaj $\alpha=\alpha_1+\alpha_2$; se \^ci tiu sumo estas ekstere la intervalo $(-\pi,\,\pi\,]$, ni devos adicii a\u u subtrahi~$2\pi$, por havi la konvenciitan $\alpha$.}
\ppr{We also remember that if $\Omega=\Omega_1\Omega_2\,$, that is, $\rho\,{\rm e}^{i\,\alpha}=(\rho_1\,{\rm e}^{i\,\alpha_1})(\rho_2\,{\rm e}^{i\,\alpha_2})$, then $\rho=\rho_1\rho_2$ and $\alpha=\alpha_1+\alpha_2$; if this sum is out the interval $(-\pi,\,\pi\,]$, we must add or subtract~$2\pi$, to have the $\alpha$ agreed.}

\ppsection[0ex]{\label{secSC}Schwarz-Christoffel}{Schwarz-Christoffel}

\ppln{Gravan familion, SC, de analitikaj mapoj $w(z)$ studis Schwarz kaj Christoffel, sendepende. Por priskribi \^gin, ni komencas kun difino: funkcio $w(z)$ estas dirita SC se \^gia {\it deriva\^\j o} havas formon \mbox{\cite[p. 550]{Hildebrand}}}
\pprn{Schwarz and Christoffel, independently, \mbox{studied} an important family (SC) of analytical maps $w(z)$. To describe it, we start with a definition: a function $w(z)$ is said SC if its {\it derivative} has the form \mbox{\cite[p. 550]{Hildebrand}}}

\bea                                                                                \label{dwdz}
\frac{\dd w}{\dd z}= \frac{C}{(z-x_1)^{k_1}(z-x_2)^{k_2}\cdots(z-x_n)^{k_n}}\,, 
\eea

\ppln{kie $C$ estas kompleksa a\u u reala konstanto, kaj $x_i$ kaj $k_i$ estas realaj konstantoj. Notu ke $\dd w/\dd z$ estas analitika, kaj ke punktoj $x_i$, en reala akso de ebeno $z$, estas singularaj punktoj por la funkcio (\ref{dwdz}); pli specife, ili kutime estos bran\^c-punktoj.}
\pprn{where $C$ is a constant complex or real, and the $x_i$ and $k_i$ are real constants. Note that $\dd w/\dd z$ is analytical, and that the points $x_i$, in the real axis of plane $z$, are singular points for the function (\ref{dwdz}); more specifically, they will generally be branch points.}

\ppl{La funkcio $\dd w/\dd z$ determinas en ebeno $w$ la {\it formon} (angulojn kaj distancojn) de la imago de ajn figuro en ebeno $z$. La funkcio indikas anka\u u la {\it orienti\^gon} de la imago, sed ne indikas la lokon de la imago en ebeno $w$. Vidu figuron~\ref{160505a}.}
\ppr{The function $\dd w/\dd z$ determines, in the plane $w$, the {\it form} (angles and distances) of the image of any figure in the plane $z$. The function also indicates the {\it orientation} of the image, but does not indicate the localization of the image in the plane $w$. See figure~\ref{160505a}.}

\begin{figure}[H]
\centerline{\epsfig{file=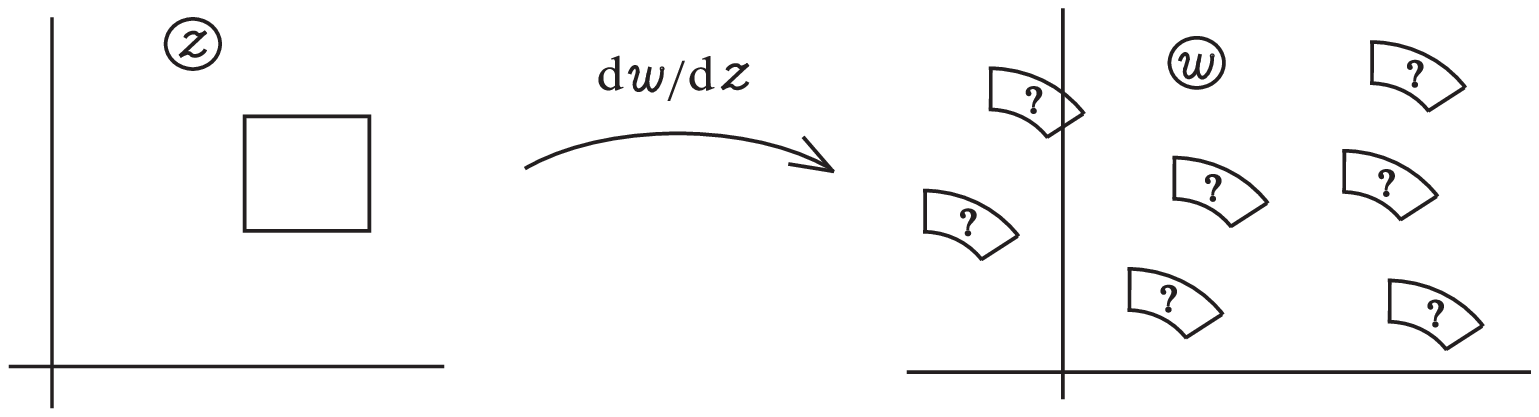,width=12cm}} 
\caption{Tutsola, la deriva\^{\j}o $\dd w/\dd z$ ne decidas kiun, el la figuroj en ebeno $w$, estas la imago de la kvadrato en ebeno $z$, per la mapo $w(z)$.
\newline
Figure~\ref{160505a}: Alone, the derivative $\dd w/\dd z$ does not decide which, from the figures in the plane $w$, is the image of the square in the plane $z$, via the map $w(z)$.}    
\label{160505a} 
\end{figure}

\ppl{La mapo $w(z)$ estas havita per malderivo de (\ref{dwdz}), kaj adicio de nova kompleksa konstanto, $K$:}
\ppr{The map $w(z)$ is obtained by indefinite integration of (\ref{dwdz}), and addition of a new complex constant, $K$:}
\bea                                                                             \label{intdwdz}
w(z)=C\int \frac{\dd z}{(z-x_1)^{k_1}(z-x_2)^{k_2}\cdots(z-x_n)^{k_n}}+K\,. 
\eea 

\ppln{Ekv.~(\ref{intdwdz}) elektas, el la nefinia kvanto de eblaj solvoj de (\ref{dwdz}), tion kion ni volas.}
\pprn{The eq. (\ref{intdwdz}) selects, from the infinity of possible solutions of (\ref{dwdz}), the one that we want.}

\ppl{En figuro~\ref{perceba}, la punkto $x_i$ estas en la akso $x$; notu, ke se la punkto $z$ estas en la supera duon-ebeno de $z$, tio estas, $\Im(z)>0$, tial \^ciam okazos 
$0<\protect\angle(z-x_i)<\pi$; kaj notu, ke se $z$ estas en la suba duon-ebeno, t.e., $\Im(z)<0$, tial \^ciam okazos $-\pi<\protect\angle(z-x_i)<0$. Tiuj asertoj veri\^gas se nur la konvencio uzita estas $\alpha\in(-\pi,\,\pi\,]$.}
\ppr{In figure~\ref{perceba}, the point $x_i$ is in the $x$-axis; note that if a point $z$ is in the upper half-plane of $z$, that is, $\Im(z)>0$, it always happens \mbox{$0<\protect\angle(z-x_i)<\pi$}; and note that if $z$ is in the lower half-plane, that is, $\Im(z)<0$, it always happens  $-\pi<\protect\angle(z-x_i)<0$. These assertions are true only if the convention used is $\alpha\in(-\pi,\,\pi\,]$.} 

\begin{figure}[H]
\centerline{\epsfig{file=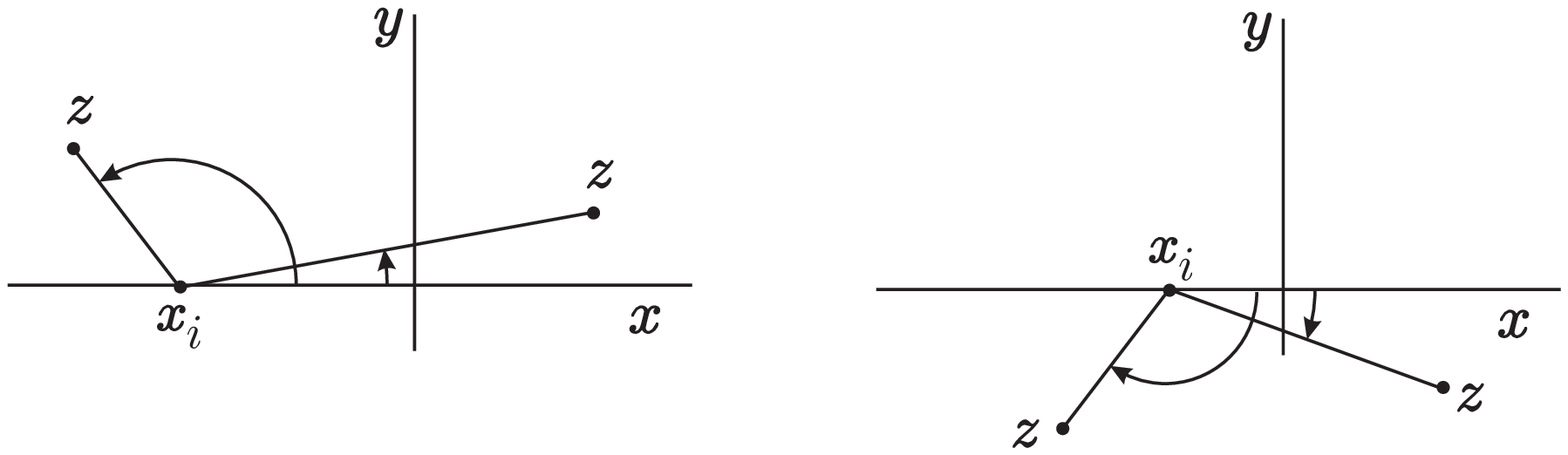,width=10cm}} 
\caption{Maldekstre, $z$ estas en la supera duon-ebeno, okazante $0<$\protect\angle$(z-x_i)<\pi$; dekstre, $z$ estas en la suba duon-ebeno, okazante $-\pi<$\protect\angle$(z-x_i)<0$.
\newline
Figure~\ref{perceba}: On the left, $z$ is in the upper half-plane, causing $0<$\protect\angle$(z-x_i)<\pi$; on the right, $z$ is in the lower half-plane, causing $-\pi<$ \protect\angle $(z-x_i)<0$.} 
\label{perceba} 
\end{figure}

\ppl{Ni konvencias, ke la $x_i$ estu aran\^gitaj tiel ke $x_1<x_2<\cdots<x_n$\,, kaj ni limigas nian studon al okazoj $-2\leq k_1+k_2+\cdots+k_n\leq2$\,.}
\ppr{We agree that the $x_i$ be ordinated so that $x_1\!<\!x_2\!<\!\cdots\!<\!x_n$\,, and we limit our study to the cases $-2\leq k_1+k_2+\cdots+k_n\leq2$\,.}

\ppl{En la du sekvantaj sekcioj ni montros, ke \^ciu SC mapas orientitan horizontalan rekton $R^+$, ku\^santa iomete {\it super} linio $y=0$, en orientita linio farita de rektaj segmentoj, en la ebeno $w$. Plue, tiu sama esprimo $w(z)$ mapas orientitan horizontalan rekton $R^-$, ku\^santa iomete {\it sub} linio $y=0$, en {\it alia} linio anka\u u farita de rektaj segmentoj. Ni opinias, ke la analizo de la mapoj de $R^+$ kaj de $R^-$ estas la kulmino de la studo de transformoj de Schwarz-Christoffel.}
\ppr{In the next two sections we show that every SC maps an oriented horizontal straight line $R^+$, laying slightly {\it above} the line $y=0$, into an oriented line made of straight segments, in the plane $w$. More, that same expression $w(z)$ maps an oriented horizontal straight line $R^-$, laying slightly {\it below} the line $y=0$, into {\it another} oriented line also made of straight segments. In our opinion, the analysis of the maps of $R^+$ and of $R^-$ is the highest point in the study of Schwarz-Christoffel transformations.}

\ppl{Indas rememori, ke `\^ciuj punktoj' en ne-finio en kompleksaj ebenoj estas unu sola punkto, skribita $\infty$. En nia teksto, la imago de tiu punkto $z=\infty$ estos skribita $w_\infty$. En la mapo de la supera duon-ebeno, la imago de najbaro de punkto $x_i$ skribi\^gos $w_i$, kontra\u ue en la mapo de suba duon-ebeno \^gi skribi\^gos $W_i$.} 
\ppr{It is worth remember that `all the points' at infinity in complex planes are one single point, written $\infty$. In our text, the image of this point $z=\infty$ will be written $w_\infty$. In the map of the upper half-plane, the image of the neighborhood of a point $x_i$ will be written $w_i$, while in the map of the lower half-plane it will be written $W_i$.}

\ppsection[0ex]{\label{superior}Mapo SC de rekto $R^+$}{Map of straight line $R^+$}

\ppln{Ni nun priskribos la mapon de la supera duon-ebeno de $z$, farita per elektita funkcio SC. Tiu duon-ebeno enhavas nur punktojn kun $y>0$, do \^gi ne enhavas la realan akson $y=0$. Ni unue ser\^cas la formon de la imago, per $w(z)$, de horizontala rekto $R^+$ ku\^santa {\it infinitezime super} la rekto $y=0$ de la ebeno $z$. Por tio, ni supozas punkton $z$ iranta en la rekto $z=x+i\,\epsilon$, estante $\epsilon$ {\it pozitiva} infinitezimo, de $x=-\infty$ al $x=+\infty$, kiel figuro~\ref{caminho} montras. Kaj ni konstruos la imagon $w(R^+)$ de tiu movado, en la ebeno $w$.}
\pprn{We shall now describe the map of the upper half-plane of $z$, made by a selected function SC. This half-plane contains only points with $y>0$, so it does not contain the real axis $x$. We first search the form of the image, under $w(z)$, of a horizontal straight line $R^+$ laying {\it infinitesimally above} the line $y=0$ of the plane $z$. To that end, we suppose a point $z$ going along the straight line $z=x+i\,\epsilon$, with $\epsilon$ a {\it positive} infinitesimal, from $x=-\infty$ to $x=+\infty$, as figure~\ref{caminho} shows. And we shall construct the image $w(R^+)$ of that motion, in the plane $w$.}

\begin{figure}[H]  
\centerline{\epsfig{file=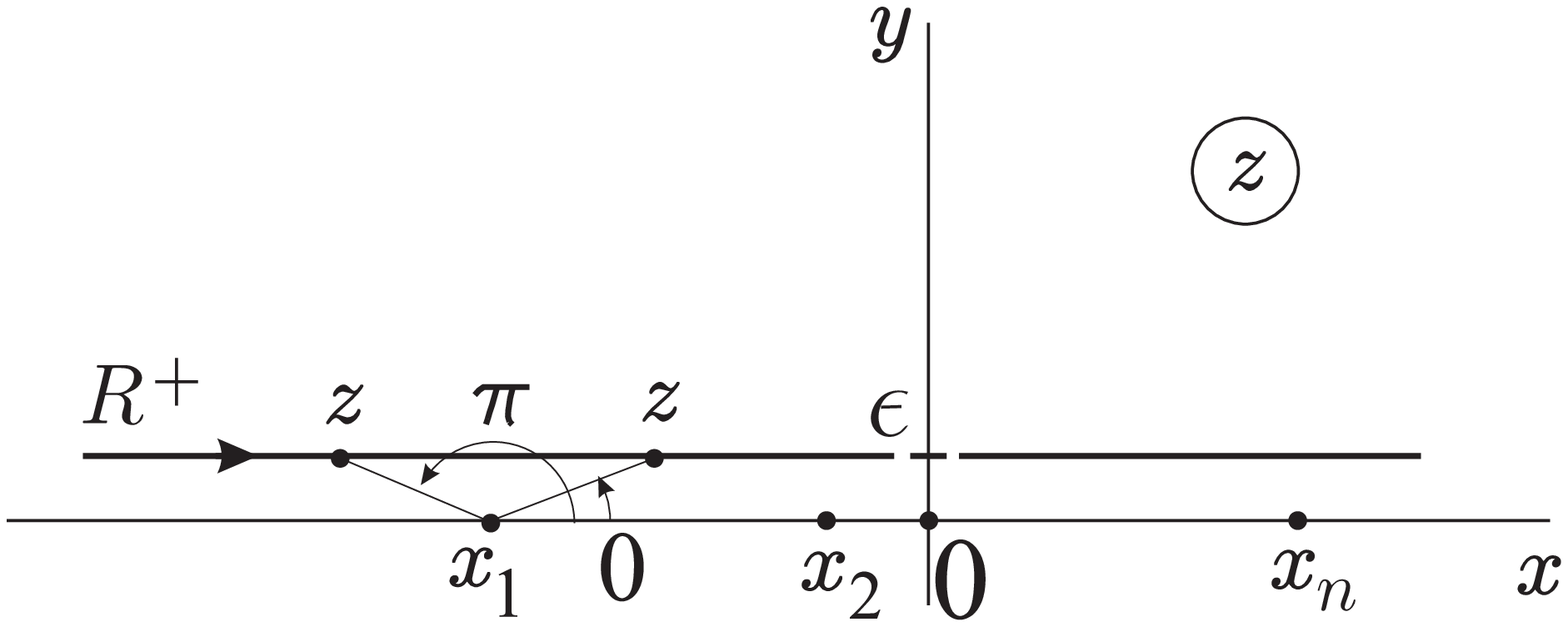,width=7cm}\hspace{4em}} 
\caption{Punkto $z$ trakuras rekton $R^+$ en la signalita direkto, estante $\epsilon$ {\it pozitiva} infinitezimo; kiam $z$ pasas super $x_1$, la angulo de $(z-x_1)$ malmilde \^san\^gas de $\pi$ al $0$.
\newline
Figure~\ref{caminho}: Point $z$ goes along the straight line $R^+$ in the indicated direction, with $\epsilon$ {\it positive} infinitesimal; when $z$ overpasses $x_1$, the angle of $(z-x_1)$ changes abruptly from $\pi$ to $0$.} 
\label{caminho} 
\end{figure}

\ppl{Komence ni konsideras la angulojn de amba\u u flankoj de la deriva\^\j o (\ref{dwdz}), elektante (tio gravas!) la valoron $k_i$\protect\angle$(z-x_i)$ por la angulo \protect\angle $(z-x_i)^{k_i}$\,:}
\ppr{We initially consider the angles of both sides of the derivative (\ref{dwdz}), choosing (this is important!) the value $k_i$\protect\angle$(z-x_i)$ for the angle \protect\angle $(z-x_i)^{k_i}$\,:}
\bea                                                                                 \label{ang}
\protect\angle[\dd w(z)/\dd z]=\protect\angle C -k_1\protect\angle(z-x_1) -k_2\protect\angle(z-x_2) -\cdots -k_n\protect\angle(z-x_n)\,.
\eea

\ppln{Dum $z=x+i\,\epsilon$ trakuras de $x=-\infty$ al tre proksima al $x_1$, \^ciuj anguloj $\protect\angle(z-x_i)$ en (\ref{ang}) valoras $\pi$, do}
\pprn{While $z=x+i\,\epsilon$ goes from \mbox{$x=-\infty$} to very near $x_1$, all angles \mbox{$\protect\angle(z-x_i)$} in (\ref{ang}) are $\pi$, so}
\bea                                                                                \label{ang1}
\protect\angle[\dd w(z)/\dd z]=\protect\angle C -\pi\,\Sigma k_i \hspace{5mm} {\rm se}\;\;x<x_1\,; 
\eea

\ppln{tial \protect\angle$[\dd w(z)/\dd z]$ estas konstanta en la intervalo $x\in(-\infty,\,x_1)$. Tio implicas, ke la imaga kurbo $w(z)$ estu rekta segmento en la komenca peco $w_\infty\rightarrow w_1$, kun {\it orienti\^go} $\alpha_0=\protect\angle C-\pi\,\Sigma k_i$. La apendico~\ref{der} pravigas tiun valoron por la orienti\^go.} 
\pprn{thus \protect\angle$[\dd w(z)/\dd z]$ is constant in the interval $x\in(-\infty,\,x_1)$. This implies that the image line $w(z)$ be a straight segment in the initial part $w_\infty\rightarrow w_1$, with {\it orientation} $\alpha_0=\protect\angle C-\pi\,\Sigma k_i$. The appendix~\ref{der} justifies this value for the orientation.}

\ppl{Figuroj~\ref{caminho} kaj \ref{161117a} montras, ke kiam $z$ pasas super $x_1$, la angulo \protect\angle$(z-x_1)$ bruske malkreskas de $\pi$ al nulo, kontra\u ue la aliaj anguloj \protect\angle$(z-x_i)$ ne \^san\^gas; tio okazigas bruskan \^san\^gon $k_1\pi$ de la direkto de $w(R^+)$ en punkto $w_1$. Notu, ke tiu \^san\^go estas pozitiva se $k_1>0$, kaj estas negativa se $k_1<0$.}
\ppr{Figures~\ref{caminho} and \ref{161117a} show that, when $z$ overpasses $x_1$, the angle \protect\angle$(z-x_1)$ abruptly decreases from $\pi$ to zero, while the other angles \protect\angle$(z-x_i)$ do not change; this causes an abrupt change $k_1\pi$ of the direction of the image $w(R^+)$ at point $w_1$. Note that this change is positive if $k_1>0$, and is negative if $k_1<0$.}

\ppl{Post $z$ pasas super $x_1$ kaj iras en la direkto de $x_2$, \^ciuj anguloj \mbox{$\protect\angle(z-x_i)$} estas konstantaj; tio farigas, ke la imago $w_1\rightarrow w_2$ denove estu rekta segmento.}
\ppr{After $z$ overpasses $x_1$ and proceeds in the direction to $x_2$, all angles \mbox{$\protect\angle(z-x_i)$} are constant; this makes the image $w_1\rightarrow w_2$ again be a straight segment.}

\ppl{Resume, la paso de $z$ super iu $x_i$ faras la orienti\^gon de la nova rekta segmento en ebeno $w$ \^san\^gi $k_i\pi$ relative al la orienti\^go de la anta\u ua segmento \mbox{\cite[p. 552]{Hildebrand}.}}
\ppr{Summing up, the passing of $z$ over any $x_i$ makes the orientation of the new straight segment in the plane $w$ change $k_i\pi$ relatively to the orientation of the previous segment \mbox{\cite[p. 552]{Hildebrand}.}}

\ppl{Fine, post $z$ pasas super $x_n$, la lasta $x_i$ de la sinsekvo, okazas ke \^ciuj anguloj \mbox{$\protect\angle(z-x_i)$} estas nulaj, do}
\ppr{At the end, after $z$ overpasses $x_n$, the last $x_i$ of the sequence, all angles \mbox{$\protect\angle(z-x_i)$} are null, so}
\bea                                                                                \label{angn}
\protect\angle[\dd w(z)/\dd z]=\protect\angle C  \hspace{5mm} {\rm se}\;\;x>x_n\,; 
\eea

\ppln{tio signifas, ke la orienti\^go de la imaga kurbo $w(z)$ de la rekto $R^+$ en la fina peco $w_n\rightarrow w_\infty$ estas $\alpha_n=\protect\angle C$. Vidu figuron~\ref{161117a}.}
\pprn{this means that the orientation of the image line $w(z)$ of the straight line $R^+$ in the final part $w_n\rightarrow w_\infty$ is $\alpha_n=\protect\angle C$. See figure~\ref{161117a}.}

\begin{figure}[H]
\centerline{\epsfig{file=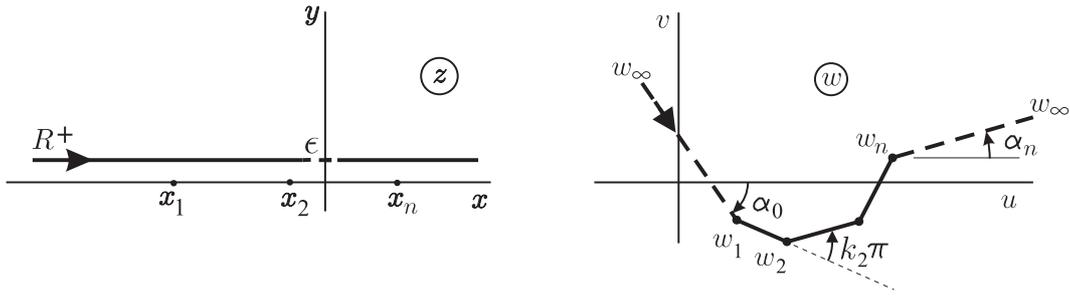,width=14cm}} 
\caption{Orientita rekto $R^+$ ku\^sas iomete super la rekto $y=0$\,, kiu entenas $n$ punktojn $x_i$. En la ebeno $w$, oni havas la imagon de $R^+$ per mapo (\ref{intdwdz}). En \^ciu vertico $w_i$ la \^san\^go de direkto estas $k_i\pi$\,. La anguloj \hspace{1mm} $\alpha_0=$ \protect\angle $C-\pi\,\Sigma k_i$ \hspace{1mm} kaj \hspace{1mm} $\alpha_n=$ \protect\angle $C$ \hspace{1mm} indikas la orienti\^gon de la komenca kaj de la fina peco de la imago.
\newline 
Figure~\ref{161117a}: The oriented straight line $R^+$ lays just above the straight line $y=0$\,, which contains $n$ points $x_i$. In the plane $w$, one has the image of $R^+$ via a map (\ref{intdwdz}). At each vertex $w_i$ the change of direction is $k_i\pi$\,. The angles \hspace{1mm} $\alpha_0=$ \protect\angle $C-\pi\,\Sigma k_i$ \hspace{1mm} and \hspace{1mm} $\alpha_n=$ \protect\angle $C$ \hspace{1mm} indicate the orientation of the initial part and of the final part of the image.} 
\label{161117a} 
\end{figure}

\ppl{Havante \^ciujn \^san\^gojn $k_i\pi$ de orienti\^go de tiu imaga kurbo, ni konsideru 4 malsamajn tipojn de mapo $w(R^+)$, la\u u la valoro $|\Sigma\,k_i\pi|$ estu plieta ol $\pi$, egala al $\pi$, inter $\pi$ kaj $2\pi$, a\u u egala al $2\pi$:}
\ppr{Having all changes $k_i\pi$ of orientation of that image curve, let's consider 4 different types of map $w(R^+)$, according the value $|\Sigma\,k_i\pi|$ is less than $\pi$, equal to $\pi$, between $\pi$ and $2\pi$, or equal to $2\pi$:}

\ppln{\noindent$\bullet$ tipo {\bf a}: se $|\Sigma\,k_i|<1$\,; 

\noindent$\bullet$ tipo {\bf b}: se $|\Sigma\,k_i|=1$\,; 

\noindent$\bullet$ tipo {\bf c}: se $1<|\Sigma\,k_i|<2$\,; kaj 

\noindent$\bullet$ tipo {\bf d}: se $|\Sigma\,k_i|=2$\,. 
\newline Figuro~\ref{contornos} montras unu imagon por \^ciu tipo, kaj sekcioj de \ref{ttipoa} al \ref{ttipod} donos unu ekzemplon de \^ciu tipo.}
\pprn{\noindent$\bullet$ type {\bf a}: if $|\Sigma\,k_i|<1$\,; 

\noindent$\bullet$ type {\bf b}: if $|\Sigma\,k_i|=1$\,; 

\noindent$\bullet$ type {\bf c}: if $1<|\Sigma\,k_i|<2$\,; and 

\noindent$\bullet$ type {\bf d}: if $|\Sigma\,k_i|=2$\,. 
\newline Figure~\ref{contornos} shows one image for each type, and sections from \ref{ttipoa} to \ref{ttipod} will give one example of each type.}

\begin{figure}[H]
\centerline{\epsfig{file=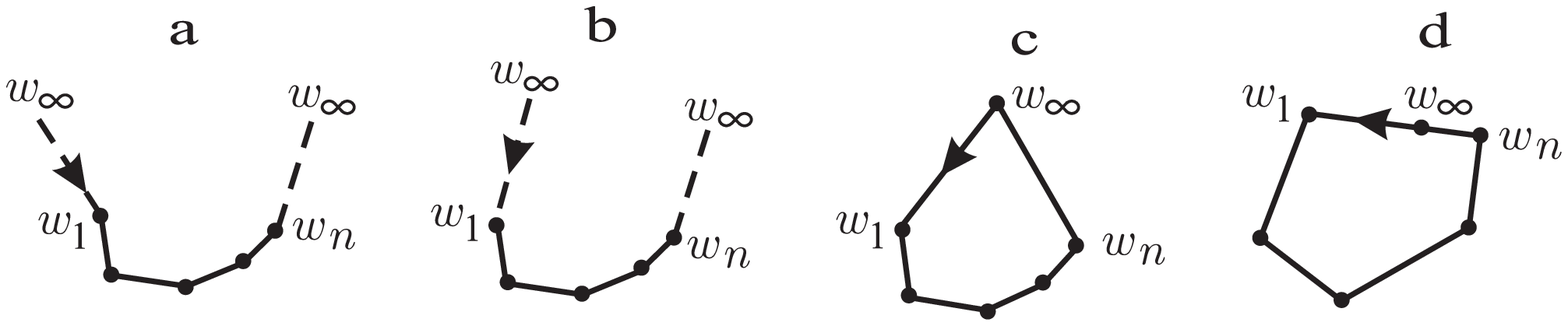,width=13cm}} 
\caption{Tipoj de imago de orientita rekto $R^+$. En tipoj {\bf a} kaj {\bf b} la imago punkto $w_\infty$ estas en ne-finio; tio implicas, ke la pecoj $w_\infty\rightarrow w_1$ kaj $w_n\rightarrow w_\infty$ estas ne-finie longaj. En tipoj {\bf c} kaj {\bf d} la imago $w_\infty$ estas en finia regiono de ebeno $w$.
\newline 
Figure~\ref{contornos}: Types of image of the oriented straight line $R^+$. In types {\bf a} and {\bf b} the image point $w_\infty$ is at infinity; this implies the parts $w_\infty\rightarrow w_1$ and $w_n\rightarrow w_\infty$ be infinitely long. In types {\bf c} and {\bf d} the image point $w_\infty$ is in the finite region of plane $w$.}    
\label{contornos} 
\end{figure}

\ppl{Por kompletigi la imagon de la rekto $R^+$, ni nun kalkulas la longojn de la rektaj segmentoj, uzante (\ref{intdwdz}),}
\ppr{To complete the image of the straight line $R^+$, we now calculate the lengths of the straight segments, using (\ref{intdwdz}),}

\bea                                                                                 \label{mod}
|w_{i+1}-w_i|=|C|\left|\int_{x_i}^{x_{i+1}}\frac{\dd x}{(x-x_1)^{k_1}(x-x_2)^{k_2}\cdots(x-x_n)^{k_n}}\right|.  
\eea

\ppl{\^Car la funkcio $w(z)$ estas kontinua, la imagoj de la aliaj rektoj $y=$\,konst\,$>0$ najbaraj al rekto $y=0$ estas similaj al la imago de rekto $R^+$, kiel ni konstatos en sekcioj de \ref{ttipoa} al \ref{ttipod}.}
\ppr{Since the function $w(z)$ is continuous, the images of the other straight lines $y=$\,const\,$>0$ neighbor to the straight line $y=0$ are similar to the image of $R^+$, as we shall verify in sections from \ref{ttipoa} to \ref{ttipod}.}

\ppsection[0.6ex]{\label{inferior}Mapo SC de rekto $R^-$}{Map of straight line $R^-$}

\ppln{Nun ni havigu la mapon de suba duon-ebeno de $z$ uzante la saman esprimon (\ref{intdwdz}) uzita por la supera duon-ebeno. Simile kiel anta\u ue, ni komence ser\^cas la imagon de horizontala rekto $R^-$ ku\^santa je infinitezime {\it negativa} ordinato $-\epsilon$ de ebeno $z$. Vidu figuron~\ref{porbaixo}.}
\pprn{Let's now get the map of the lower half-plane of $z$ using the same expression (\ref{intdwdz}) used for the upper half-plane. Similarly as before, we initially search the image of a horizontal straight line $R^-$ laying in {\it negative} infinitesimal ordinate $-\epsilon$ of plane $z$. See figure~\ref{porbaixo}.}

\begin{figure}[H]
\centerline{\epsfig{file=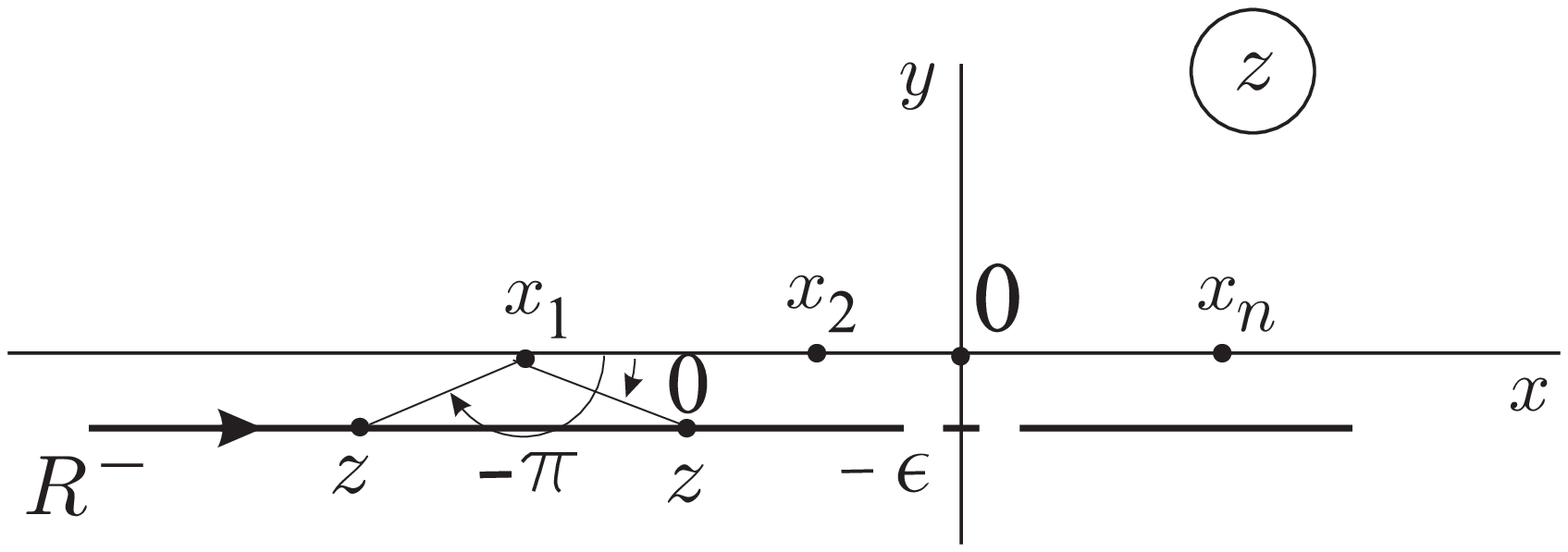,width=8cm}\hspace{4.9em}} 
\caption{Kiam $z$ pasas sub $x_1$ en la signalita direkto, estante $-\epsilon$ {\it negativa} infinitezimo, la angulo de $(z-x_1)$ bruske \^san\^gas de $-\pi$ al $0$.
\newline 
Figure~\ref{porbaixo}: When $z$ passes under $x_1$ in the indicated direction, being $-\epsilon$ a {\it negative} infinitesimal, the angle of $(z-x_1)$ abruptly changes from $-\pi$ to $0$.}     
\label{porbaixo} 
\end{figure}

\ppl{Ni nun supozas punkton $z=x-i\,\epsilon$ iranta dekstren en rekto $R^-$, ekde $x=-\infty$. Inter $x=-\infty$ kaj $x=x_1$ \^ciuj anguloj \protect\angle$(z-x_i)$ valoras $-\pi$, do la imago de komenca peco havas la konstantan orienti\^gon (\ref{ang})}
\ppr{We now suppose a point $z=x-i\,\epsilon$ going rightwards in the straight line $R^-$, since $x=-\infty$. Between $x=-\infty$ and $x=x_1$ all angles \mbox{\protect\angle$(z-x_i)$} are $-\pi$, so the image of the initial part of $R^-$ has orientation (\ref{ang}) constant}
\bea                                                                                \label{ang2}
\protect\angle[\dd w(z)/\dd z]=\protect\angle C +\pi\,\Sigma\,k_i \hspace{5mm} {\rm se}\;\;x<x_1\,. 
\eea

\ppl{Kiam $z$ pasas sub $x_1$, la angulo de $(z-x_1)$ bruske kreskas de $-\pi$ al nulo, kontra\u ue la anguloj de la aliaj $(z-x_i)$ ne \^san\^gas; tio implicas \^san\^gon $-k_1\pi$ de orienti\^go de $w(R^-)$ en punkto $W_1$. Ni notas, ke tiu \^san\^go havas saman absolutan valoron, sed kontra\u uan signumon, al la \^san\^go en la responda imaga punkto $w_1$ de rekto~$R^+$.}
\ppr{When $z$ passes under $x_1$, the angle of \mbox{$(z-x_1)$} abruply increases from $-\pi$ to zero, while the angles of the other $(z-x_i)$ do not change; this implies a change $-k_1\pi$ in the orientation of $w(R^-)$ at point $W_1$. We note that this change has same absolute value, but opposite sign, as that of the change in the corresponding image point $w_1$ of straight line  $R^+$.}

\ppl{Da\u urigante ni konstatas, ke en \^ciu imaga punkto $W_i$ de rekto $R^-$ la \^san\^go de orienti\^go havas saman absolutan valoron, sed kontra\u uan signumon, al de la responda imaga punkto $w_i$ de rekto~$R^+$. Post $z$ pasas sub $x_n$, \^ciuj anguloj \mbox{\protect\angle$(z-x_i)$} estas nulaj, do (\ref{ang}) reskribi\^gas}
\ppr{Proceeding, we see that at every image point $W_i$ of the straight line $R^-$ the change of direction has same absolute value, but opposite sign, as that on the corresponding image point $w_i$ of straight line $R^+$. After $z$ passes under $x_n$, all angles \mbox{$\protect\angle(z-x_i)$} are null, so (\ref{ang}) rewrites}

\vspace{-3mm}
\bea                                                                               \label{angn2}
\protect\angle[\dd w(z)/\dd z]=\protect\angle C  \hspace{5mm} {\rm if}\;\;x>x_n\,; 
\eea

\ppln{tial ni vidas, ke la orienti\^go de la fina peco $W_n\rightarrow w_\infty$ de la imaga kurbo de rekto $R^-$ valoras $\protect\angle C$. Tiu valoro koincidas kun la orienti\^go (\ref{angn}) de la fina peco $w_n\rightarrow w_\infty$ de la imago de rekto~$R^+$.}
\pprn{thus we see that the orientation of the final part $W_n\rightarrow w_\infty$ of the image of straight line $R^-$ is $\protect\angle C$. This value coincides with the orientation (\ref{angn}) of the final part $w_n\rightarrow w_\infty$ of the image of straight line $R^+$.}

\ppl{Per analizo de (\ref{mod}) ni konstatas, ke \^ciu longo $\;|W_{i+1}-W_i|\;$ egalas la respondan longon \mbox{$|w_{i+1}-w_i|$}. Do, la formo de imago de rekto $R^-$ havi\^gas per iu reflekto de la imago de rekto $R^+$. \^Car la orienti\^goj de la imagoj de peco  $x_n\rightarrow\infty$ en rektoj $R^+$ kaj $R^-$ koincidas (amba\u u estas $\alpha_n=\protect\angle C$), tial la reflekto okazas paralele al fina rekta peco $w_n\rightarrow w_\infty$\,. Vidu figuron~\ref{reflex2}.}
\ppr{By analysis of (\ref{mod}) we see that each size \mbox{$|W_{i+1}-W_i|$} equals the corresponding size $|w_{i+1}-w_i|$. Thus, the form of the image of straight line $R^-$ is obtained via some reflection of the image of straight line $R^+$. Since the orientations of the images of the part $x_n\rightarrow\infty$ in straight lines $R^+$ and $R^-$ coincide (both are $\alpha_n=\protect\angle C$), the reflection occurs parallel to the final straight part $w_n\rightarrow w_\infty$\,. See figure~\ref{reflex2}.}

\begin{figure}[H]
\centerline{\epsfig{file=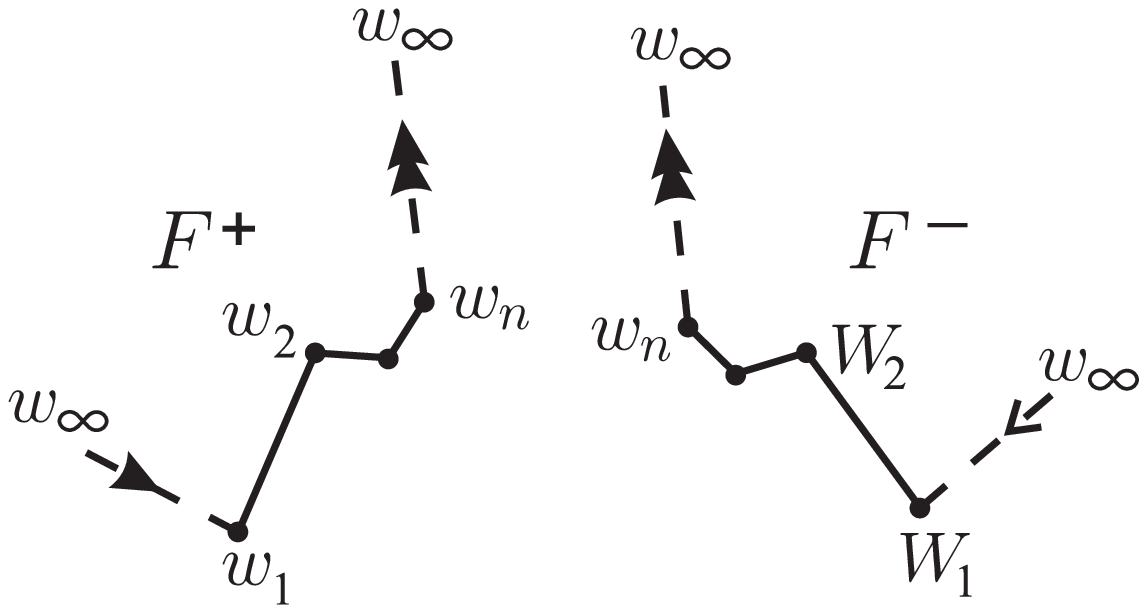,width=6cm}} 
\caption{Se $F^+$ estas la {\it formo} de la imago de rekto $R^+$, tial $F^-$ estos la {\it formo} de la imago de responda rekto $R^-$\,; oni vidas, ke $F^-$ estas reflekto de $F^+$ en la direkto $w_n\rightarrow w_\infty$\,.
\newline
Figure~\ref{reflex2}: If $F^+$ is the {\it form} of the image of straight line $R^+$, then $F^-$ is the {\it form} of the corresponding straight line $R^-$\,; one sees that $F^-$ is refletion of $F^+$ in the direction $w_n\rightarrow w_\infty$\,.} 
\label{reflex2} 
\end{figure}

\ppl{La sekvantaj kvar sekcioj prezentas ekzemplojn de la kvar tipoj de transformo $w(z)$ priskribitaj en sekcio~\ref{superior}. Sekcio~\ref{cinem} prezentas praktikan aplikon.}
\ppr{The next four sections give examples of the four types of transformation $w(z)$ described in section~\ref{superior}. Section~\ref{cinem} presents a practical application.}

\ppsection[0.1ex]{\label{ttipoa}Ekzemplo de tipo a}{An example of type a}

\ppln{Se en (\ref{dwdz}) ni elektas $n=1$\,, $C=1$\,, $x_1=0$\,, $k_1=1/2$\,, ni havos}
\pprn{If in (\ref{dwdz}) we choose $n=1$\,, $C=1$\,, \mbox{$x_1=0$\,,} $k_1=1/2$\,, we shall have}

\vspace{-3mm}
\bea                                                                               \label{tipoa}
\frac{\dd w}{\dd z}=\frac{1}{z^{1/2}}\,, \hspace{3mm} \Rightarrow \hspace{3mm} w(z)=2\,z^{1/2}+K\,; 
\eea

\ppln{elektante $K=0$, ni havigas imagojn $w(z)$ kiel en figuro~\ref{160426e}\,.}
\pprn{choosing $K=0$, we obtain images $w(z)$ as in figure~\ref{160426e}\,. }

\begin{figure}[H]
\centerline{\epsfig{file=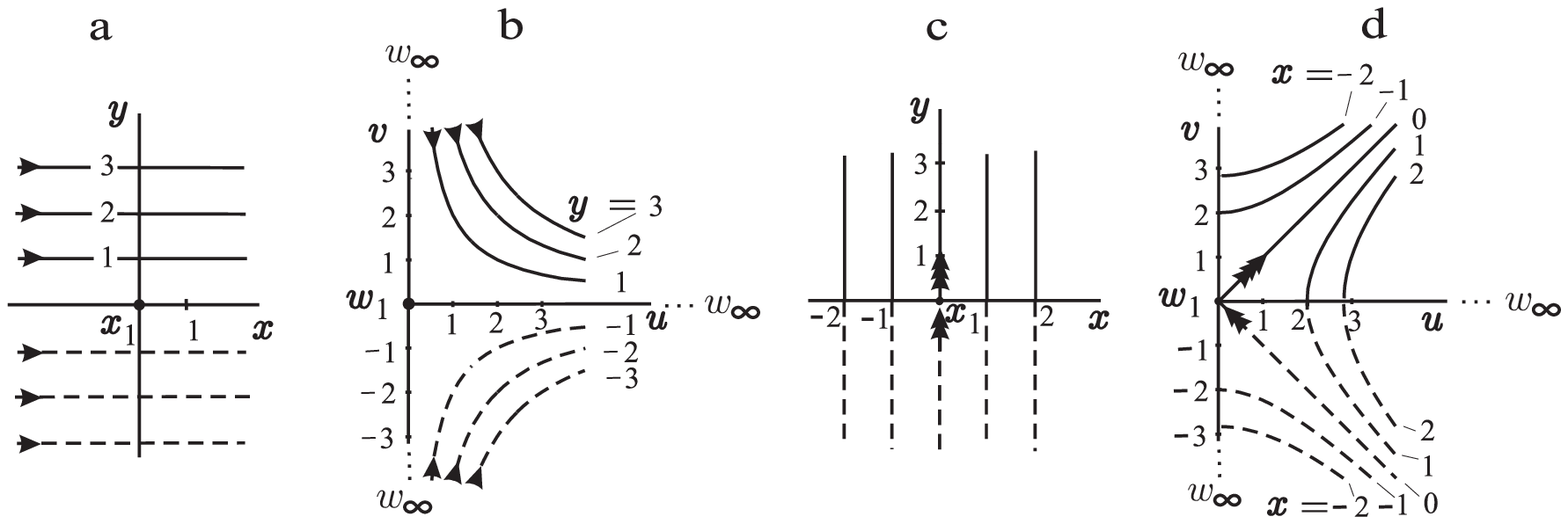,width=\textwidth}} 
\caption{{\bf a} kaj {\bf c} montras rektojn $y=\,$konst kaj $x=\,$konst en ebeno $z$\,; {\bf b} kaj {\bf d} montras la respondajn imagojn en ebeno $w$ per la transformo $w(z)=2\,z^{1/2}$\,. 
\newline
Figure~\ref{160426e}: {\bf a} and {\bf c} show straight lines $y=\,$const and $x=\,$const in plane $z$\,; {\bf b} and {\bf d} show the corresponding images in plane $w$ via the transformation $w(z)=2\,z^{1/2}$\,.}
\label{160426e} 
\end{figure}

\ppl{En figuro~\ref{160426e}{\bf b} ni vidas, ke la imago de $R^+$ trakuras de $v=+\infty$ al origino, poste iras al $u=+\infty$; ni vidas anka\u u, ke la imago de $R^-$ trakuras de $v=-\infty$ al origino, poste iras al $u=+\infty$. Notu, en figuro~\ref{160426e}{\bf d}, ke imagoj de rektoj $x=\,{\rm konst\,}<0$ estas malkontinuaj, kaj notu ke la imago de rekto $x=0$ estas zigzaga en punkto $w_1$\,; kaj fine notu ke imagoj de rektoj $x=\,{\rm konst}\,>0$ estas kontinuaj kaj ne-zigzagaj. Kontra\u ue, la imagoj de \^ciuj horizontalaj rektoj ($y={\rm konst}$) estas kontinuaj; tio estos tre interesa kiam ni uzos tiajn mapojn por priskribi fluadon de likvoj, kiel en sekvantaj sekciojn.}
\ppr{In figure~\ref{160426e}{\bf b} we see that the image of $R^+$ runs from $v=+\infty$ till the origin, then goes to $u=+\infty$; we also see that the image of $R^-$ runs from $v=-\infty$ till the origin, then goes to $u=+\infty$. Note, in figure~\ref{160426e}{\bf d}, that the images of the lines $x=\,{\rm const\,}<0$ are discontinuous, and note that the image of the line $x=0$ brakes at point $w_1$; further note that the images of the lines $x=\,{\rm const}\,>0$ are continuous and do not brake. Oppositely, the images of all horizontal lines ($y={\rm const}$) are continuous; this will be very interesting when we shall use such maps to describe flow of liquids, as in sections to follow.}

\ppl{Eble interesas akompani, pa\^so post pa\^so, konforman mapon ekde komenca formo \^gis la fina formo. Figuro~\ref{161124b} skizas kvar imagojn, per $\dd w/\dd z=1/z^{k_1}$, de horizontalaj rektoj $y=\,$konst, por $k_1$ kreskanta ekde $0$ \^gis $1/2$.}
\ppr{It is perhaps interesting to follow, step by step, a conformal map since the initial form till the final one. Figure~\ref{161124b} sketches four images, for $\dd w/\dd z=1/z^{k_1}$, of horizontal lines $y=\,$const, for $k_1$ increasing from $0$ to $1/2$.}
\begin{figure}[H]
\centerline{\epsfig{file=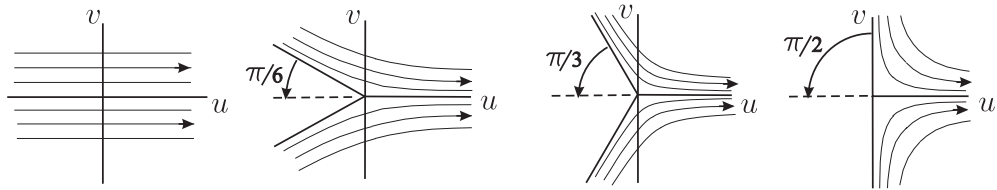,width=13cm}} 
\caption{Imagoj de rektoj $y=\,$konst per $\dd w/\dd z=1/z^{k_1}$, kun $k_1$ sinsekve $0$, $1/6$, $1/3$ kaj $1/2$.
\newline
Figure~\ref{161124b}: Images of straight lines $y=\,$const via $\dd w/\dd z=1/z^{k_1}$, with $k_1=$ $0$, $1/6$, $1/3$ and $1/2$.} 
\label{161124b} 
\end{figure}
\ppsection[0.6ex]{\label{ttipob}Ekzemplo de tipo b}{An example of type b}

\ppln{En (\ref{dwdz}), kun $n=2$\,, $C=1$\,, $x_1=-1$\,, $x_2=1$\,, $k_1=k_2=\frac{1}{2}$, ni havigas}
\pprn{In (\ref{dwdz}), with $n=2$\,, $C=1$\,, $x_1=-1$\,, $x_2=1$\,, $k_1=k_2=\frac{1}{2}$, we obtain}

\bea                                                                               \label{tipob}
\frac{\dd w}{\dd z}=\frac{1}{(z+1)^{1/2}(z-1)^{1/2}}\,,\hspace{3mm} \Rightarrow \hspace{3mm} w(z)=\cosh^{-1}(z)+K\,. 
\eea

\ppln{Ni elektas $K=0$ kaj havas imagojn kiel en figuroj~\ref{160502f} kaj \ref{160502g}.}
\pprn{We choose $K=0$ and have images as in \mbox{figures~\ref{160502f} and \ref{160502g}.}}

\begin{figure}[H]
\centerline{\epsfig{file=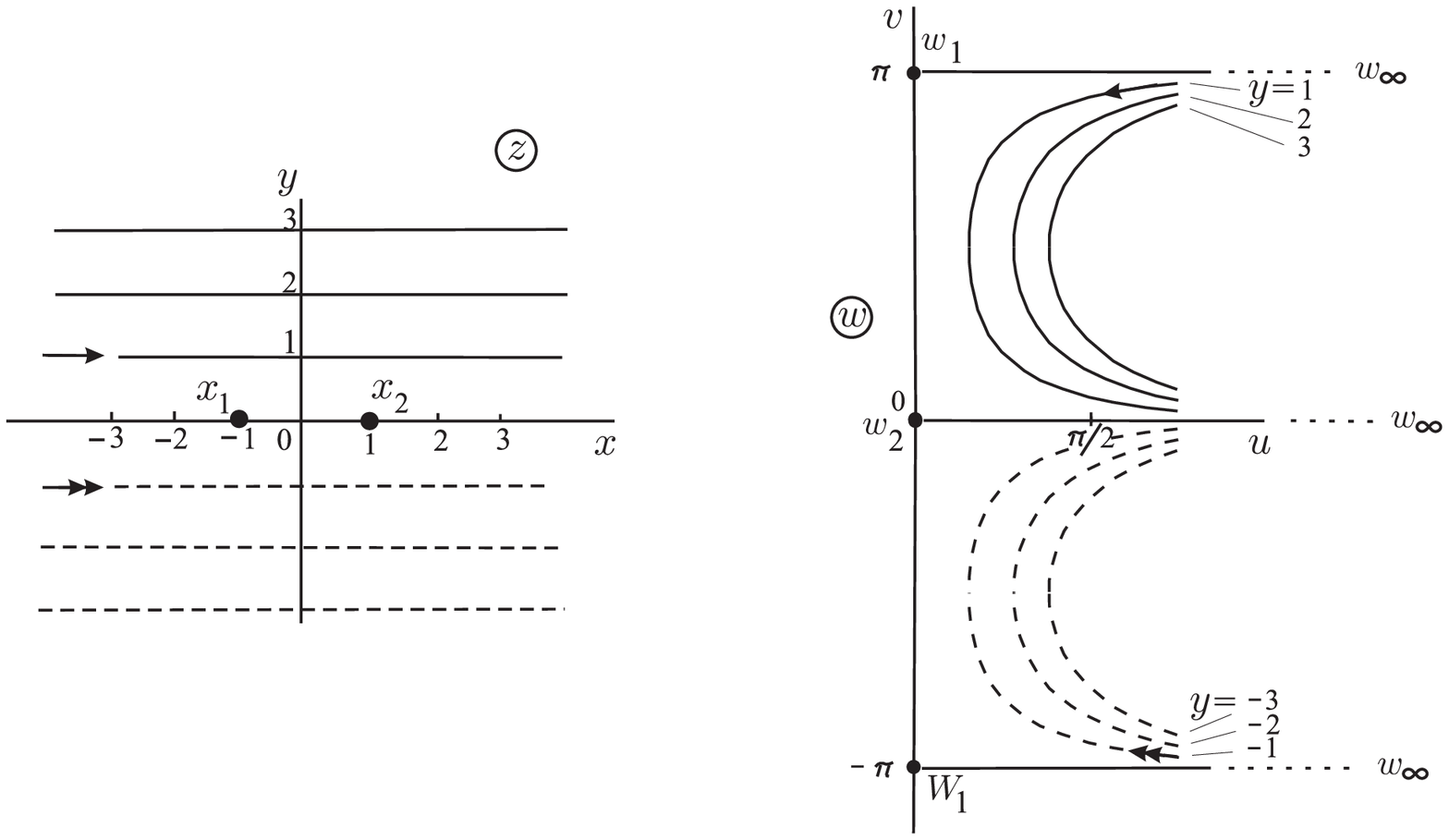,width=14cm}} 
\caption{Horizontalaj rektoj $y=\,$konst en ebeno $z$ kaj iliaj imagoj per mapo $w(z)=\cosh^{-1}(z)$.
\newline
Figure~\ref{160502f}: Straight lines $y=\,$const in plane $z$ and their images via the map $w(z)=\cosh^{-1}(z)$.}
\label{160502f} 
\end{figure}

\begin{figure}[H]
\centerline{\epsfig{file=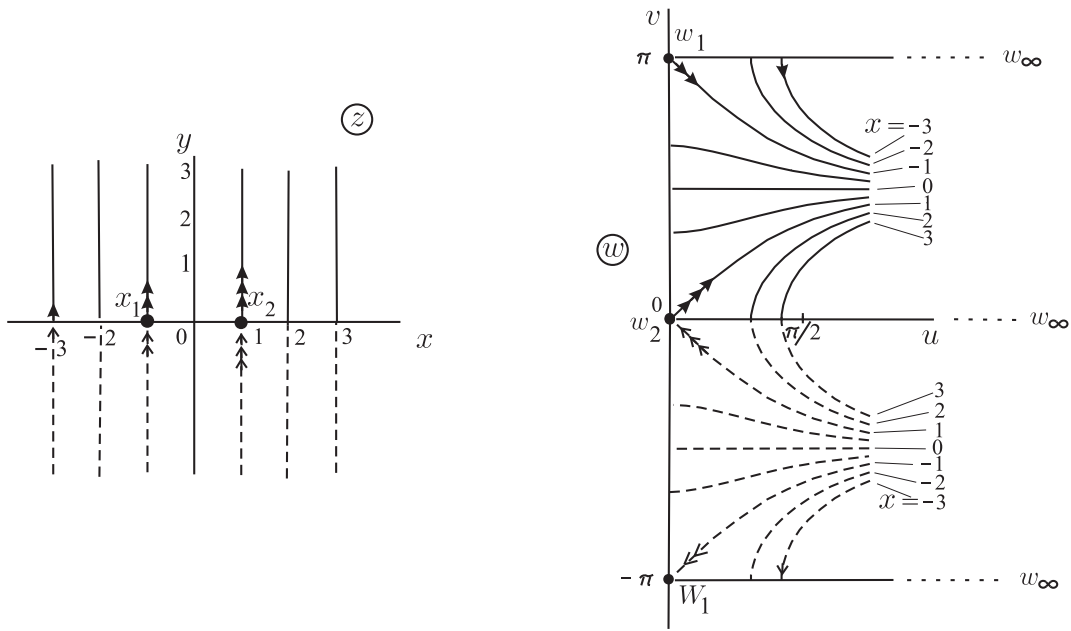,width=14cm}} 
\caption{Vertikalaj rektoj $x=\,$konst en ebeno $z$ kaj iliaj imagoj per la mapo $w(z)=\cosh^{-1}(z)$.
\newline
Figure~\ref{160502g}: Straight lines $x=\,$const in plane $z$ and their images via the map $w(z)=\cosh^{-1}(z)$.}
\label{160502g} 
\end{figure}

\ppl{En figuro~\ref{160502f} ni vidas, ke la$\,$ imago$\,$ de$\,$ $R^+\,$ venas$\,$ horizontale$\,$ de \mbox{$u=+\infty$} \^gis $u=0$, kun $v=\pi$, da\u urigas en la akso $v$ \^gis la origino, poste iras horizontale al $u=+\infty$ en la akso $u$. La imago de $R^-$ venas horizontale de $u=+\infty$ \^gis $u=0$, kun $v=-\pi$, da\u urigas en la akso $v$ \^gis la origino, kaj iras horizontale al $u=+\infty$ en la akso $u$. Tiel, la imago de $R^+$ kaj $R^-$ kun $x>x_2$ koincidas je la pozitiva parto de akso $u$.}
\ppr{In figure~\ref{160502f} we see that the image of $R^+$ comes horizontally from $u=+\infty$ till $u=0$, with $v=\pi$, stays in the axis $v$ till the origin, then goes horizontally to $u=+\infty$ in the axis $u$. The image of $R^-$ comes horizontally from $u=+\infty$ till $u=0$, with $v=-\pi$, stays in the axis $v$ till the origin, then goes horizontally to $u=+\infty$ in axis $u$. Thus the images of $R^+$ and $R^-$ with $x>x_2$ coincide in the positive part of $u$-axis.}

\ppl{Kiel en anta\u ua ekzemplo (tipo a), \^ciuj horizontaj rektoj en ebeno $z$ havas kontinuajn imagojn, kiel montras figuro~\ref{160502f}. Denove, ne \^ciuj vertikalaj rektoj en ebeno $z$ havas kontinuajn imagojn, kiel montras figuro~\ref{160502g}; nur imagoj de rektoj kun $x\geq1$ estas kontinuaj.}
\ppr{Similarly as in the previous example (type a), all horizontal lines in plane $z$ have continuous images, as figure~\ref{160502f} shows. And again, not all vertical lines in plane $z$ have continuous image, as figure~\ref{160502g} shows; only the images of vertical lines with $x\geq1$ are continuous.}

\ppsection[0.6ex]{\label{ttipoc}Ekzemplo de tipo c}{An example of type c}

\ppln{Se en (\ref{dwdz}) ni elektas $n=2$\,, $C=1$\,, $x_1=-1$\,, $x_2=1$\,, $k_1=k_2=\frac{2}{3}$\,, ni havos}
\pprn{If in (\ref{dwdz}) we choose $n=2$\,, $C=1$\,, $x_1=-1$\,, $x_2=1$\,, $k_1=k_2=\frac{2}{3}$\,, we shall have}

\bea                                                                               \label{tipoc}
\frac{\dd w}{\dd z}=\frac{1}{(z+1)^{2/3}(z-1)^{2/3}} \hspace{3mm} \Rightarrow \hspace{3mm} w(z)=-{\rm e}^{i\pi/3}\,z\,F\left(1/2,\:2/3\,; 3/2\,;\,z^2\right)+K\,,  
\eea 

\ppln{kie $F$ estas la hipergeometria funkcio donata per la serio de Taylor \mbox{\cite[p. 179]{Hildebrand}}}
\pprn{where $F$ is the hypergeometric function given by the Taylor series \mbox{\cite[p. 179]{Hildebrand}}}

\bea                                                                               \label{hiper}
F(a,b\,;c\,;\zeta)=1+\frac{a\,b}{c}\,\zeta + \frac{a(a+1)\,b(b+1)}{c(c+1)}\,\frac{\zeta^2}{2!} + \frac{a(a+1)(a+2)\,b(b+1)(b+2)}{c(c+1)(c+2)}\,\frac{\zeta^3}{3!}+\cdots\,. 
\eea

\ppln{Plue, elektante la konstanto $K=\alpha\,{\rm e}^{i\,2\pi/3}$ kun $\alpha=F\left(\,\frac{1}{2}\:,\,\frac{2}{3}\:;\,\frac{3}{2}\:;\,1\right)\approx2,1$
ni havigas, post iun da kalkuloj, $w_2=-w_\infty$ kaj la mapojn montritajn en figuroj~\ref{140904f} kaj \ref{140904g}.}
\pprn{Further choosing constant $K=\alpha\,{\rm e}^{i\,2\pi/3}$ with $\alpha=F\left(\,\frac{1}{2}\:,\,\frac{2}{3}\:;\,\frac{3}{2}\:;\,1\right)\approx2,1$
we obtain, after some calculus, $w_2=-w_\infty$ and the maps shown in \mbox{figures~\ref{140904f} and \ref{140904g}.}}

\begin{figure}[H]
\centerline{\epsfig{file=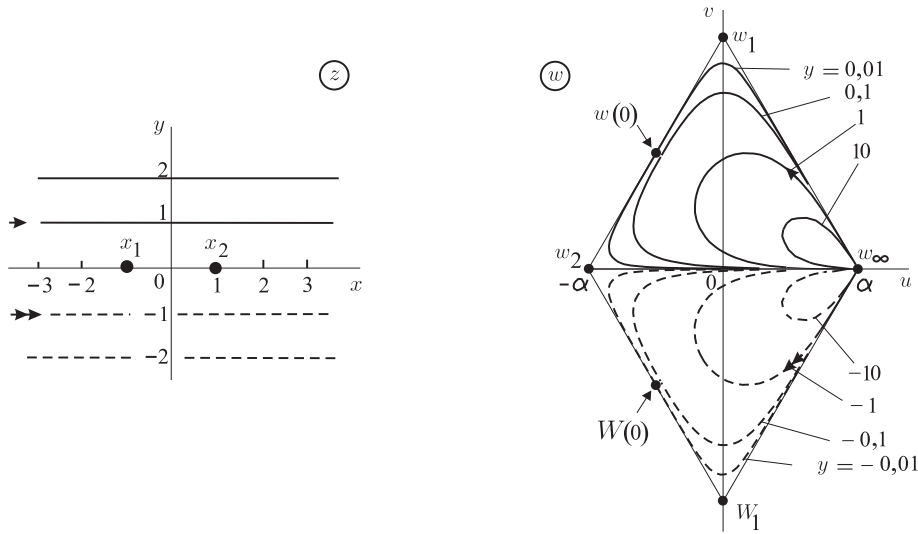,width=12cm}} 
\caption{Horizontalaj rektoj $y=\,$konst en ebeno $z$ kaj iliaj imagoj en ebeno $w$ per mapo~(\ref{tipoc}).
\newline
Figure~\ref{140904f}: Straight lines $y=\,$const in plane $z$ and their images in plane $w$ via map~(\ref{tipoc}).}
\label{140904f} 
\end{figure}

\begin{figure}[H]
\centerline{\epsfig{file=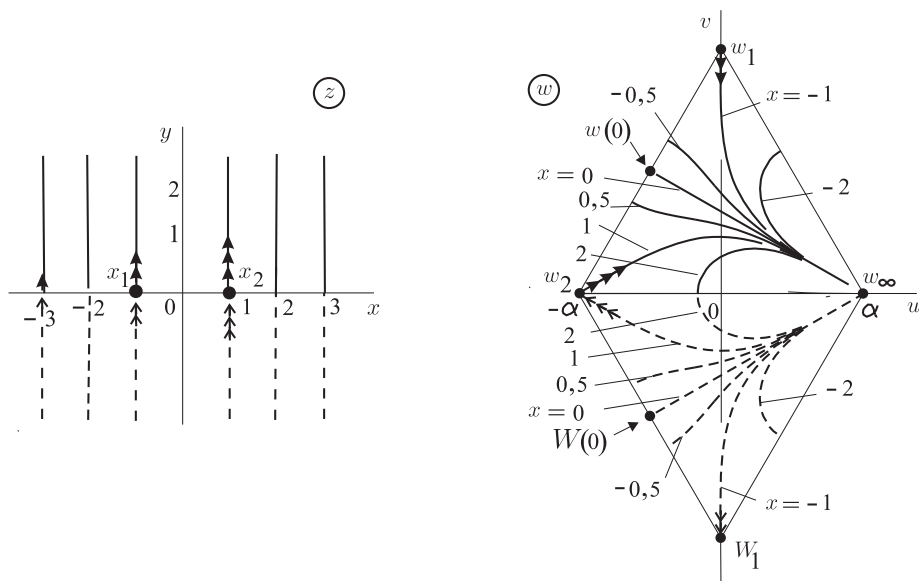,width=12cm}} 
\caption{Vertikalaj rektoj $x=\,$konst en ebeno $z$ kaj iliaj imagoj en ebeno $w$ per mapo~(\ref{tipoc}).
\newline
Figure~\ref{140904g}: Vertical straight lines $x=\,$const in plane $z$ and their images in plane $w$ via map~(\ref{tipoc}).}
\label{140904g} 
\end{figure}

\ppl{En figuro~\ref{140904f} ni vidas, ke la imagoj de $R^+$ kaj de $R^-$ estas du egallateraj trianguloj. La horizontalaj lateroj koincidas, kaj estas la imagoj de $R^+$ kaj $R^-$ kun $x>x_2$.}
\ppr{In figure~\ref{140904f} we see that the images of $R^+$ and $R^-$ are two equilateral triangles. The horizontal sides coincide, and are the images of $R^+$ and $R^-$ with $x>x_2$.}

\ppl{Kiel en anta\u uaj du ekzemploj, \^ciu horizontala rekto en ebeno $z$ havas kontinuan imagon, kiel montras figuro~\ref{140904f}. Denove, ne \^ciu vertikala rekto havas kontinuan imagon, kiel montras figuro~\ref{140904g}.}
\ppr{Similarly as in the previous two examples, every horizontal straight line in plane $z$ has continuous image, as figure~\ref{140904f} shows. And again, not all vertical straight lines in plane $z$ have continuous image, as figure~\ref{140904g} shows.}

\ppsection[0.6ex]{\label{ttipod}Ekzemplo de tipo d}{An example of type d}

\ppln{Se en (\ref{dwdz}) kun $n=3$ ni uzos $x_1=-1$\,, $x_2=0$\,, $x_3=1$\,, $k_1=\frac{3}{4}$\,, $k_2=\frac{1}{2}$\,, kaj $k_3=\frac{3}{4}$\,, ni havos}
\pprn{If in (\ref{dwdz}) with $n=3$ we use $x_1=-1$\,, $x_2=0$\,, $x_3=1$\,, $k_1=\frac{3}{4}$\,, $k_2=\frac{1}{2}$\,, and $k_3=\frac{3}{4}$\,, we shall have}
\bea                                                                               \label{tipod}
\frac{\dd w}{\dd z}=\frac{C}{(z+1)^{3/4}z^{1/2}(z-1)^{3/4}}\hspace{2mm} \Rightarrow w(z)=-2\,C\,{\rm e}^{i\pi/4}\,z^{1/2}\,F\left(\,1/4,\: 3/4\,; \,5/4\,; \,z^2\right)+K\,, 
\eea
\ppln{estante $K$ konstanto, kaj estante $F$ la hipergeometria funkcio kiel en (\ref{hiper}). Ni elektas $C=1$ kaj $K=i\,\sqrt{2}\,F(\,\frac{1}{4},\;\frac{3}{4}\,;\frac{5}{4}\,;\,1)\approx2,62\,i$, kaj havigas imagojn kiel en figuroj~\ref{160511e} kaj \ref{160511f}.}
\pprn{being $K$ a constant and being $F$ the hypergeometric function as in (\ref{hiper}). We choose $C=1$ and $K=i\,\sqrt{2}\,F(\,\frac{1}{4},\;\frac{3}{4}\,;\frac{5}{4}\,;\,1)\approx2,62\,i$, and obtain images as in figures~\ref{160511e} and \ref{160511f}.}  

\begin{figure}[H]
\centerline{\epsfig{file=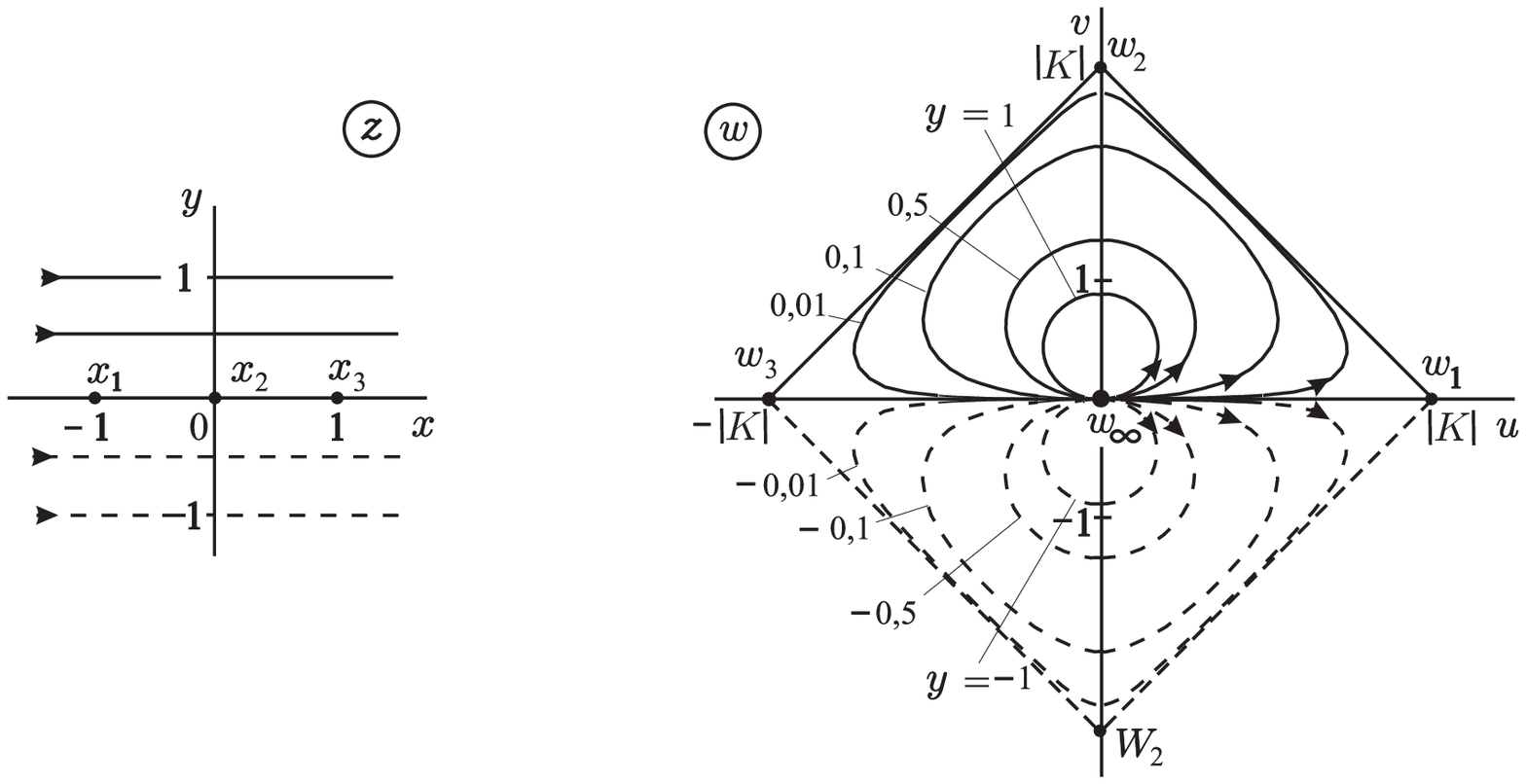,width=14cm}} 
\caption{Horizontalaj rektoj $y=\,$konst en ebeno $z$, kaj iliaj imagoj per la mapo~(\ref{tipod}). 
\newline
Figure~\ref{160511e}: Horizontal straight lines $y=\,$const in plane $z$, and their images via the map~(\ref{tipod}).}  
\label{160511e} 
\end{figure}

\begin{figure}[H]
\centerline{\epsfig{file=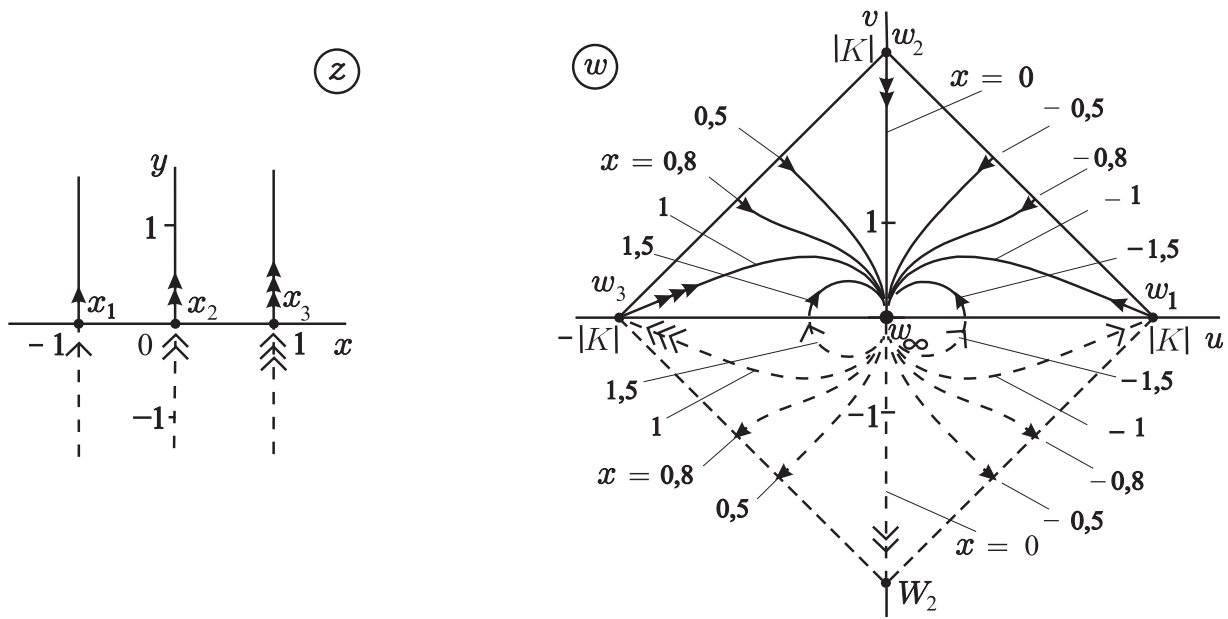,width=16cm}} 
\caption{Vertikalaj rektoj $x=\,$konst en ebeno $z$, kaj iliaj imagoj per la mapo~(\ref{tipod}) 
\newline
Figure~\ref{160511f}: Vertical straight lines $x=\,$const in plane $z$, and their images via the map~(\ref{tipod}).}  
\label{160511f} 
\end{figure}

\ppl{En figuro~\ref{160511e} vidu, ke la imagoj de $R^+$ kaj $R^-$ estas ortaj trianguloj. La horizontalaj hipotenuzoj koincidas, kaj estas la imagoj de $R^+$ kaj $R^-$ kun $x<x_1$ kaj $x>x_3$. }
\ppr{In figure~\ref{160511e}, see that the images of $R^+$ and $R^-$ are rectangular triangles. The horizontal hypotenuses coincide, and are the images of $R^+$ and $R^-$ with $x<x_1$ and $x>x_3$.}

\ppl{Kiel en anta\u uaj tri  ekzemploj, \^ciu horizontala rekto en ebeno $z$ havas kontinuan imagon, kiel montras figuro~\ref{160511e}. Denove, ne \^ciu vertikala rekto en ebeno $z$ havas kontinuan imagon, kiel montras figuro~\ref{160511f}.}
\ppr{Similarly as in the previous three examples, every horizontal line in plane $z$ has continuous image, as figure~\ref{160511e} shows. And again, not all vertical lines in plane $z$ have continuous image, as figure~\ref{160511f} shows.}

\ppsection[0.6ex]{\label{cinem}Piliero en rivero}{A pillar in a river}

\ppln{En figuroj de \ref{160426e} al \ref{160511e} ni vidas, ke la imago, en ebeno $w$, de linioj $y={\rm konst}$ el ebeno $z$, simili\^gas al linioj de fluo de iu fluido. Tiu simileco ne estas akcidenta: vere, baza literaturo \mbox{\cite[p. 313]{Hildebrand}}
montras ke SC mapo generas liniojn de fluo de fluido limigita per rektaj segmentoj. Ni nun prezentas ekzemplon de mapo SC, kiu simulas la movadon de akvo en kanalo kun ortangula piliero en \^gia bordo; la dimensioj de tiu piliero estas malgrandaj kompare al la lar\^go de la kanalo. Ni perceptos, en figuro~\ref{160128b}, ke \^ci tiu ekzemplo ta\u ugas anka\u u por piliero en mezo de kanalo.}
\pprn{In figures from \ref{160426e} to \ref{160511e} we see that the image of the horizontal lines $y={\rm const}$, in plane $w$, resemble lines of flow of some fluid. This resemblance is not accidental: indeed, basic texts \mbox{\cite[p. 313]{Hildebrand}} show that an SC map generates lines of flow of fluid bounded by straight segments. We now present an example of SC map which simulates the motion of water in a channel with a rectangular pillar on its boundary; the dimensions of this pillar are small in comparison with the width of the channel. We shall perceive, in figure~\ref{160128b}, that this example serves also for a pillar in the middle of a channel.}

\begin{figure}[H]
\centerline{\epsfig{file=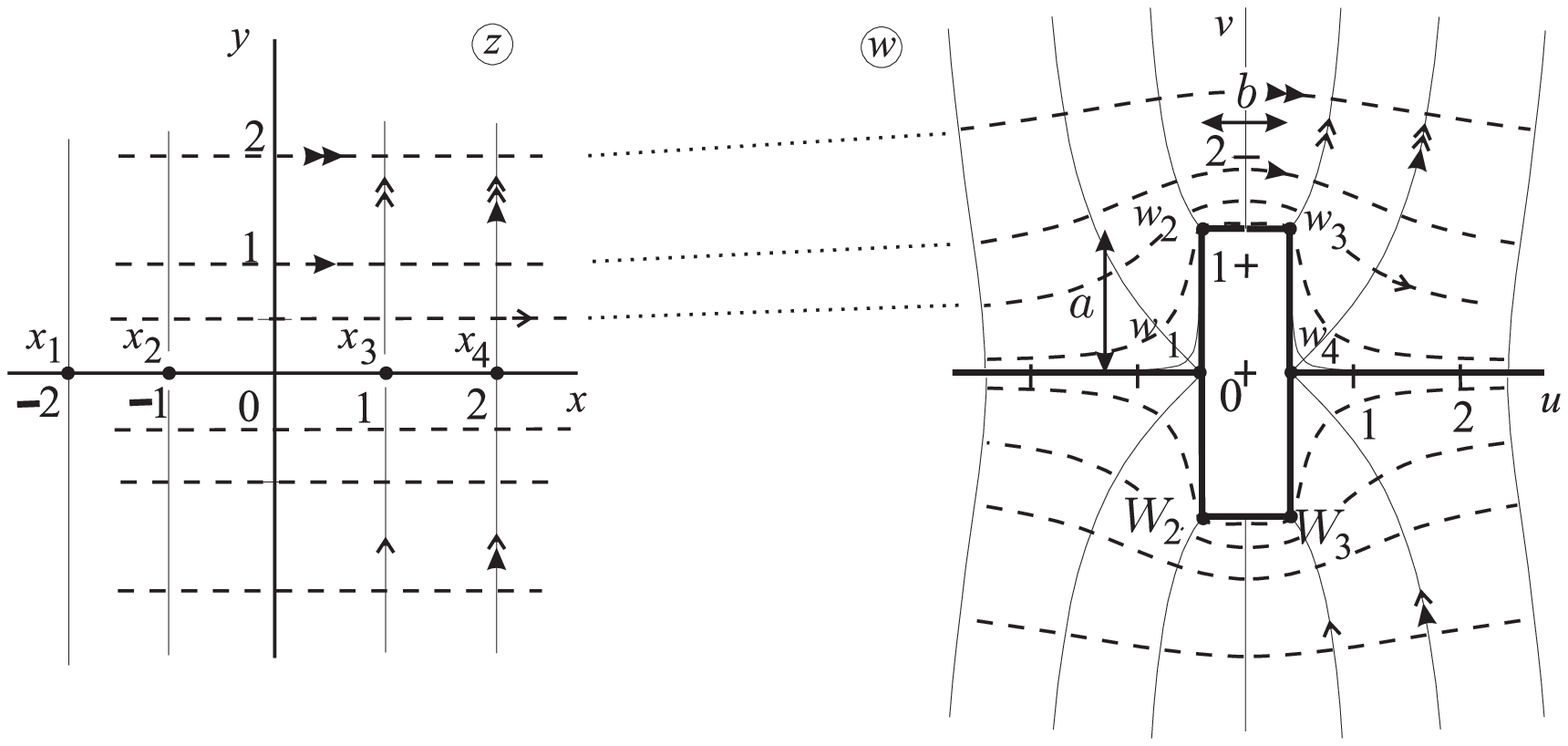,width=14cm}} 
\caption{La supera (malsupera) duon-ebeno de $z$ estas mapita en la supera (malsupera) duon-ebeno de $w$ per (\ref{murinho})\,.
\newline
Figure~\ref{160128b}: The upper (lower) half-plane of $z$ is mapped into the upper (lower) half-plane of $w$ via (\ref{murinho}).} 
\label{160128b} 
\end{figure}

\ppl{Por havi la ortan formon ni uzas kvar ortangulojn: $k_1=k_4=1/2$ kaj $k_2=k_3=-1/2$. Ni elektas konvenajn $x_i$ kaj $C$: $x_1=-2$\,, $x_2=-1$\,, $x_3=1$\,, $x_4=2$\,, $C=1$. Tiel (\ref{dwdz}) donas}
\ppr{To have the rectangular form we use four right angles: $k_1=k_4=1/2$ and $k_2=k_3=-1/2$. We choose convenient $x_i$ and $C$: $x_1=-2$\,, $x_2=-1$\,, $x_3=1$\,, $x_4=2$\,, $C=1$. So (\ref{dwdz}) gives}

\bea                                                                             \label{murinho}
\frac{\dd w}{\dd z}=\frac{1}{(z+2)^{1/2}(z+1)^{-1/2}(z-1)^{-1/2}(z-2)^{1/2}}=
\left(\frac{z^2-1}{z^2-4}\right)^{1/2}
\Rightarrow \hspace{1mm} w(z)=E(z/2,\,2)+K, \mbox{~~} 
\eea

\ppln{kie $E$ estas la nekompleta eliptika integrala\^\j o de dua tipo \mbox{\cite[p. 170]{Dwight},}}
\pprn{where $E$ is the incomplete elliptical integral of second type \mbox{\cite[p. 170]{Dwight},}}

\bea                                                                           \label{elipticas}
E(Z,\,k)=\int_0^Z\frac{\sqrt{1-k^2\zeta^2}\;\dd\zeta}{\sqrt{1-\zeta^2}}\,. 
\eea

\ppl{La grandoj $a:=v_2-v_1$ kaj $b:=u_3-u_2$ de la piliero estas}  
\ppr{The dimensions $a:=v_2-v_1$ and $b:=u_3-u_2$ of the pillar are}  

\bea                                                                               \label{hkayl}
a=\int_1^2\sqrt{\frac{x^2-1}{4-x^2}}\,\dd x\approx1,3 \hspace{4mm} {\rm and} \hspace{3mm}  b=\int_{-1}^1\sqrt{\frac{1-x^2}{4-x^2}}\,\dd x\approx0,8\,. 
\eea

\ppl{Figuro~\ref{160128b} respondas al valoro $K=1,3\,i$ en ekv.~(\ref{murinho}), kaj montras imagojn de kelkaj horizontalaj kaj vertikalaj rektoj en ebeno $z$. Tiu figuro estas pripensebla je du eblecoj: piliero je mezo de kanalo, kaj piliero ^ce bordo de kanalo. Fakte, se ni forigas la suban duon-ebenon de $w$, tiel la supera duon-ebeno montras pilieron \^ce bordo de kanalo.}
\ppr{Figure~\ref{160128b} corresponds to the value $K=1,3\,i$ in eq.~(\ref{murinho}), and shows images of some horizontal and vertical lines in plane $z$. This figure can be thought under two possibilities: pillar in middle of the channel, and pillar on the border of the channel. Indeed, if we eliminate the lower half-plane of $w$, then the upper half-plane shows a pillar on the border of the channel.}

\ppsection[-0.4ex]{\label{komentoj}Komentoj}{Comments}

\ppln{La uzado de ne-entjeraj potencoj $k_i$ en (\ref{dwdz}) kaj (\ref{intdwdz}) necesigas konvencion por angulo de komplekso. \^Ci tie ni proponis, ke anguloj estu montritaj en la intervalo $(-\pi,\,\pi]$; tio multe faciligas havigi imagojn de rektoj $R^-$ (sekcio~\ref{inferior}). Tamen, \^ci tiu konvencio ne estas universala; aliaj intervaloj ofte uzitaj estas $[0,\,2\pi)$ kaj $(-\pi/2,\,3\pi/2]$\,, ekzemple. Se iu el tiuj alternativaj konvencioj estas uzita, nia formulo (\ref{ang2}) de angulo de deriva\^{\j}o necesigas modifon, same kiel la analizo farita en sekcio~\ref{inferior}.}
\pprn{The use of non-integer powers $k_i$ in (\ref{dwdz}) and (\ref{intdwdz}) forces a convention for angle of a complex. Here we proposed that the angles be shown in the interval $(-\pi,\,\pi]$; this greately facilitates obtaining images of straight lines $R^-$ (section~\ref{inferior}). However, this convention is not universal; other intervals often used are $[0,\,2\pi)$ and $(-\pi/2,\,3\pi/2]$\,, for example. If any of these alternative conventions is used, our formula (\ref{ang2}) for angle of the derivative needs a modification, as well as the analysis made in section~\ref{inferior}.}

\ppl{Plue, la konvencio $\alpha\in(-\pi,\,\pi]$ ne sufi\^cas por precizigi la transformon $w(z)$ kiel en (\ref{dwdz}). Fakte, ankora\u u estas iuj elektoj kiojn ni devas fari por havi la solvon kion ni volas. Ekzemple, se iu eksponento $k_i$ estas $m/n$, kie $m$ kaj $n$ estas interprimoj (2 kaj 9, ekzemple), tiam $(z-x_i)^{m/n}$ havas $n$ malsamajn eblajn orienti\^gojn, el kiuj nur unu eble interesas al ni. En \^ci tiu teksto ni elektis ke la angulo $\protect\angle(z-x_i)^{k_i}$ estas $k_1\protect\angle(z-x_i)$.}
\ppr{More, the convention $\alpha\in(-\pi,\,\pi]$ does not suffice to disambiguate the transformation $w(z)$ as in (\ref{dwdz}). Indeed, there are still some choices that we must make to have the solution we want. For example, if some exponent $k_i$ is $m/n$, where $m$ and $n$ are coprimes (2 and 9, say), then $(z-x_i)^{m/n}$ has $n$ different possible orientations, among which only one possibly is of our interest. In this text we chose the angle $\protect\angle(z-x_i)^{k_i}$ to be $k_1\protect\angle(z-x_i)$.}

\ppsection[0.6ex]{Konkludo}{Conclusion}

\ppl{Komencante pri harmoniaj funkcioj, ni difinis SC mapojn. La interesa aspekto de tiaj mapoj estas ke iliaj bordoj estas rektaj segmentoj. Ni tiam eksploris plurajn eblecojn de segmentoj, kaj klasifikis ilin je kvar tipoj: a, b, c, d. Ni montris ekzemplon el \^ciu tipo kaj poste ni komentis, ke tiaj mapoj memorigas nin pri movado de fluido. Tiam ni pripensis pli interesan ekzemplon de piliero en rivero, faritan per SC mapo.}
\ppr{Starting with harmonic functions, we defined SC maps. The most interesting feature of such maps is that their boundaries are straight segments. We then explored several possibilities of segments, and classified them into four types: a, b, c, d. One example of each type was shown, and we commented that such maps resemble motion of fluids. We then imagined a more interesting example, that of a pillar in a river, via an SC map.}

\ppl{Krom tiuj fizikaj aspektoj de SC mapoj, ni detale pridiskutis la difinon de angulo en tiuj mapoj. Tiu diskuto enmiksi\^gis tra nia tuta artikolo. En sekcio~\ref{komentoj} ni revenis al tiu diskuto kaj aldonis e\^c pli da detalojn.}
\ppr{Besides these physical aspects of SC maps, we discussed in detail the definition of angle in these maps. This discussion permeated all our article. In section~\ref{komentoj} we came back to that discussion, and even presented more details.}

\ppl{En sekvanta artikolo ni plu eksploros modeladon de fluido per SC mapoj. Ni donos pli interesajn ekzemplojn kaj studos kinematikajn detalojn kiel rapido, ktp.}
\ppr{In a next article we shall explore more the modeling of fluids via SC maps. We shall give more interesting examples and shall study kinematical details such as velocities, etc.}

\appendix
\ppsection[0.6ex]{\label{der}Notoj pri $\dd w(z)/\dd z$}{Notes on $\dd w(z)/\dd z$} 
\ppln{Bonkonate, en studo de realaj funkcioj $f(x)$, la valoro de deriva\^{\j}o $\dd f(x)/\dd x$ en \^ciu punkto $x$ estas la valoro de trigonometria tangento de angulo inter rekto tangenta al kurbo $y=f(x)$ kaj akso $x$. Populare, $\dd y(x)/\dd x=\tan\theta$. Demando: kion oni povas diri pri la responda deriva\^{\j}o $\dd w(z)/\dd z$, en okazo de funkcioj~$w(z)$ analitikaj kompleksaj\,?}
\pprn{As is well known, in the study of real functions $f(x)$, the value of the derivative $\dd f(x)/\dd x$ in each point $x$ is the value of the trigonometric tangent of the angle between the straight line tangent to the curve $y=f(x)$ and the $x$-axis. Popularly, $\dd y(x)/\dd x=\tan\theta$. Question: what can be said about the corresponding derivative $\dd w(z)/\dd z$, in case of the complex analytical functions $w(z)$\,?}

\begin{figure}[H]
\centerline{\epsfig{file=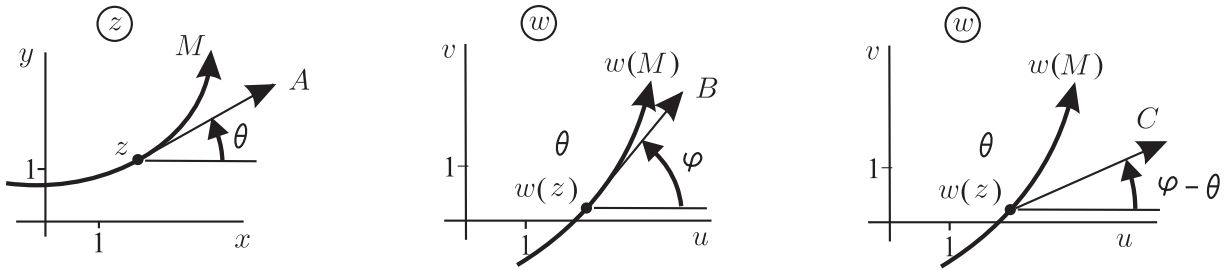,width=16cm}} 
\vspace{-1em}
\caption{Maldestre, ebeno $z$ kun orientita kurbo $M$ kaj geometria tangento $A$ en punkto $z$. Centre, ebeno $w$ kun orientita imaga kurbo $w(M)$ kaj geometria tangento $B$ en imaga punkto $w(z)$. Dekstre, la orienti\^go $C$ de la kompleksa nombro $\dd w/\dd z$ en la imaga punkto $w(z)$.  
\newline
Figure~\ref{161121c}: On the left, plane $z$ with oriented curve $M$ and the geometric tangent $A$ at point $z$. On the center, plane $w$ with the oriented image $w(M)$ and the geometric tangent $B$ at the image point $w(z)$. On the right, the orientation $C$ of the complex number $\dd w/\dd z$ at the image point $w(z)$.}
\label{161121c} 
\end{figure}

\ppl{Indas pliprecizigi la demandon. Vidu, en figuro~\ref{161121c}, orientitan kurbon $M$ en kompleksa ebeno $z$, kaj estu donata analitika funkcio $w(z)$; tiam ni havos, en la kompleksa ebeno $w$, la orientitan kurbon $w(M)$, imagon de la kurbo $M$. Elektante punkton $z$ en kurbo $M$, ni volas vida\^\j on de la kompleksa nombro $\dd w(z)/\dd z$ en la imaga punkto $w(z)$.}
\ppr{It is worth better explain the question. In figure~\ref{161121c}, see an oriented curve $M$ in complex plane $z$, and let be given an analytic function $w(z)$; then we shall have, in the complex plane $w$, the oriented curve $w(M)$, image of curve $M$. Selecting a point $z$ in curve $M$, we want a visualization of the complex number $\dd w(z)/\dd z$ at the image point $w(z)$.}

\ppl{Komence ni elektas punkton $z+\dd z$ en la kurbo $M$, najbare al la punkto $z$; la imagoj de tiuj punktoj estas $w(z)$ kaj $w(z+\dd z)$, amba\u u lokitaj en la imaga kurbo $w(M)$. Ni ser\^cas la vida\^\j on de la kvociento (kompleksa nombro) $\dd w/\dd z=[w(z+\dd z)-w(z)]/\dd z$.}
\ppr{We initially select a point $z+\dd z$ on curve $M$, neighbor to the point $z$; the images of these points are $w(z)$ and $w(z+\dd z)$, both on the image curve $w(M)$ on plane $w$. We want a visualization of the quotient (a complex number) $\dd w/\dd z=[w(z+\dd z)-w(z)]/\dd z$.}

\ppl{La infinitezima vektoro  $\dd z$ estas orientita kiel la geometria tangento $A$ al la kurbo $M$ en la punkto $z$. Simile, la infinitezima vektoro $\dd w=w(z+\dd z)-w(z)$ estas  orientita kiel la geometria tangento $B$ al la imaga kurbo $w(M)$ en la imaga punkto $w(z)$. Tial la modulo de $\dd w/\dd z$ estas la kvociento de la moduloj de tiuj infinitezimaj vektoroj, kaj la orienti\^go de $\dd w/\dd z$ estas la diferenco de la respondaj orienti\^goj; tio estas, \protect\angle$(\dd w/\dd z)=$\protect\angle$\dd w-$\protect\angle$\dd z=\varphi-\theta$. Figuro~\ref{161121c} montras la orienti\^gon $C$ de la kompleksa nombro $\dd w(z)/\dd z$.}
\ppr{The infinitesimal vector $\dd z$ has the orientation of the geometric tangent $A$ to the curve $M$ at point $z$. Similarly, the infinitesimal vector $\dd w=w(z+\dd z)-w(z)$ has the orientation of the geometric tangent $B$ to the image curve $w(M)$ at the image point $w(z)$. Thus the modulus of $\dd w/\dd z$ is the quotient of the moduli of these infinitesimal vectors, and the orientation of $\dd w/\dd z$ is the difference \mbox{between} the correspondentes orientations; that is, \protect\angle$(\dd w/\dd z)=$\protect\angle$\dd w-$\protect\angle$\dd z=\varphi-\theta$. Figure~\ref{161121c} shows the orientation $C$ of the complex number $\dd w(z)/\dd z$.}

\ppl{Do la vektoro $\dd w(z)/\dd z$ \^generale ne estas same orientita kiel la geometria tangento $B$ al imaga kurbo $w(M)$. Por ke la du orienti\^goj koincidas, bezonas okazi $\theta=0$, tio estas, $A$ bezonas esti {\it horizontala}. Okazas, ke linioj $R^+$ kaj $R^-$ en transformoj de Schwarz-Christoffel estas horizontalaj rektoj, tiam la angulo de la deriva\^{\j}o $\dd w/\dd z$ koincidas kun la klino de la imago de tiuj rektoj; \^ci tiu fakto estis uzita en sekcioj \ref{superior} kaj \ref{inferior}.}
\ppr{Thus the vector $\dd w(z)/\dd z$ has not, in general, the same orientation as the geometric tangent $B$ to the image curve $w(M)$. In order the two orientations coincide, it must happen  $\theta=0$, that is, $A$ need be {\it horizontal}. It happens that the straight lines $R^+$ and $R^-$ in Schwarz-Christoffel transformations are horizontal, so the angle of the derivative $\dd w/\dd z$ coincides with the inclination of the image of these straight lines; this fact was used in sections se\c c\~oes \ref{superior} and \ref{inferior}.}


\ppparallel{\end{Parallel}}
\end{document}